\documentclass[final,onefignum,onetabnum]{siamart220329}

\headers{Stabilization of Euler--Poincar\'e Systems with Broken Symmetry~I}{C.~Contreras and T.~Ohsawa}

\title{Controlled Lagrangians and Stabilization of Euler--Poincar\'e Mechanical Systems\\with Broken Symmetry~I: Kinetic Shaping\thanks{This paper is an expanded version of our conference paper~\cite{CoOh2021a}.
    \funding{This work was supported by NSF grant CMMI-1824798.}
  }
}

\author{
  C\'esar Contreras\thanks{Department of Mathematical Sciences, The University of Texas at Dallas, 800 W Campbell Rd, Richardson, TX 75080-3021
    (\email{cxc145430@utdallas.edu,tomoki@utdallas.edu}).
  }
  \and
  Tomoki Ohsawa\footnotemark[2]
}

\usepackage{amsfonts,amssymb}
\usepackage{verbatim}
\usepackage{graphicx}
\usepackage[caption=false]{subfig}

\crefname{subsection}{Section}{Sections}

\usepackage{enumitem}
\setlist[enumerate,1]{leftmargin=2em, label=(\roman*)}
\setlist[itemize,1]{leftmargin=1.25em}

\usepackage[group-separator={,}]{siunitx}

\usepackage{tikz}
\usetikzlibrary{arrows,snakes,backgrounds,patterns,shapes.geometric,calc}
\newcommand\centerofmass{%
    \tikz[radius=0.4em] {%
        \fill (0,0) -- ++(0.4em,0) arc [start angle=0,end angle=90] -- ++(0,-0.8em) arc [start angle=270, end angle=180];%
        \draw (0,0) circle;%
    }%
}

\newsiamremark{remark}{Remark}
\newsiamremark{example}{Example}

\usepackage[numbers,sort&compress]{natbib}

\def\defeq{\mathrel{\mathop:}=}
\def\setdef#1#2{{\left\{ #1 \ |\ #2 \right\}}}

\def\od#1#2{\dfrac{d#1}{d#2}}
\def\pd#1#2{\dfrac{\partial #1}{\partial #2}}
\def\tpd#1#2{\partial #1/\partial #2}

\def\fd#1#2{\dfrac{\delta #1}{\delta #2}}
\def\tfd#1#2{\delta #1/\delta #2}

\def\dzero#1#2{\left.\od{}{#1} #2 \right|_{#1=0}}

\def\parentheses#1{{\left(#1\right)}}
\def\brackets#1{{\left[#1\right]}}
\def\braces#1{{\left\{#1\right\}}}

\def\norm#1{{\left\|#1\right\|}}

\def\DS{\displaystyle}
\def\R{\mathbb{R}}

\def\defeq{\mathrel{\mathop:}=}
\def\eqdef{=\mathrel{\mathop:}}
\def\setdef#1#2{{\left\{ #1 \ |\ #2 \right\}}}
\def\ip#1#2{{\left\langle#1,#2\right\rangle}}

\def\diag{\operatorname{diag}}

\def\eps{\varepsilon}

\def\GL{\mathsf{GL}}
\def\SO{\mathsf{SO}}
\def\SE{\mathsf{SE}}

\def\gl{\mathfrak{gl}}

\def\so{\mathfrak{so}}
\def\se{\mathfrak{se}}

\newcommand\ad{\operatorname{ad}}

\def\PB#1#2{\left\{#1,#2\right\}}

\makeatletter
\DeclareFontFamily{OMX}{MnSymbolE}{}
\DeclareSymbolFont{MnLargeSymbols}{OMX}{MnSymbolE}{m}{n}
\SetSymbolFont{MnLargeSymbols}{bold}{OMX}{MnSymbolE}{b}{n}
\DeclareFontShape{OMX}{MnSymbolE}{m}{n}{
    <-6>  MnSymbolE5
   <6-7>  MnSymbolE6
   <7-8>  MnSymbolE7
   <8-9>  MnSymbolE8
   <9-10> MnSymbolE9
  <10-12> MnSymbolE10
  <12->   MnSymbolE12
}{}
\DeclareFontShape{OMX}{MnSymbolE}{b}{n}{
    <-6>  MnSymbolE-Bold5
   <6-7>  MnSymbolE-Bold6
   <7-8>  MnSymbolE-Bold7
   <8-9>  MnSymbolE-Bold8
   <9-10> MnSymbolE-Bold9
  <10-12> MnSymbolE-Bold10
  <12->   MnSymbolE-Bold12
}{}

\let\llangle\@undefined
\let\rrangle\@undefined
\DeclareMathDelimiter{\llangle}{\mathopen}%
                     {MnLargeSymbols}{'164}{MnLargeSymbols}{'164}
\DeclareMathDelimiter{\rrangle}{\mathclose}%
                     {MnLargeSymbols}{'171}{MnLargeSymbols}{'171}
\makeatother
\def\metric#1#2{\llangle #1, #2\rrangle}

\def\bOmega{\boldsymbol{\Omega}}
\def\bXi{\boldsymbol{\Xi}}
\def\bv{\mathbf{v}}
\def\bPi{\boldsymbol{\Pi}}
\def\bP{\mathbf{P}}
\def\bGamma{\boldsymbol{\Gamma}}

\def\bchi{\boldsymbol{\chi}}
\def\ellc{\ell_{\tau,\sigma,\rho}}

\begin{document}

\maketitle

\begin{abstract}
  We extend the method of controlled Lagrangians with kinetic shaping to those mechanical systems on semidirect product Lie groups with broken symmetry, more specifically to the Euler--Poincar\'e equations with advected parameters.
  We find a matching condition for the controlled Lagrangian for such systems whose configuration manifold is a general semidirect product Lie group $\mathsf{G} \ltimes V$.
  Our motivating examples are a bottom-heavy underwater vehicle and a top spinning on a movable base.
  Their configuration space is the special Euclidean group $\mathsf{SE}(3) = \mathsf{SO}(3) \ltimes \mathbb{R}^{3}$, where the $\mathsf{SE}(3)$-symmetry is broken by the gravity.
  The controls resulting from the matching condition stabilize unstable equilibria of these examples.
  Furthermore, the matching helps us find additional dissipative controls that asymptotically stabilize those unstable equilibria.
\end{abstract}

\begin{keywords}
  Stabilization; controlled Lagrangians; Euler--Poincar\'e mechanical systems; broken symmetry; semidirect product
\end{keywords}

\begin{AMS}
  34H15, 37J25, 70E17, 70H33, 70Q05, 93D05, 93D15
\end{AMS}

\section{Introduction}
\subsection{Motivating examples}
\begin{figure}[hbtp]
  \centering
  \subfloat[Bottom-heavy underwater vehicle]{
    \includegraphics[width=0.45\linewidth]{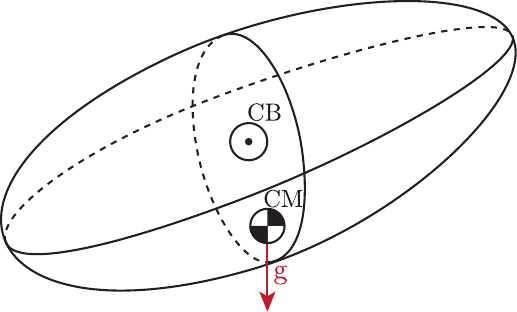}
    \label{fig:uwv}
  }
  \quad
  \subfloat[Heavy top spinning on movable base]{
    \begin{tikzpicture}[scale=.75]
      \node[draw, cylinder, shape aspect=3, rotate=90, minimum height=1cm, minimum width=3cm] (c1)  at (0,1){};
      \draw[fill] (0,1.5) circle [radius=0.05];

      \node[semitransparent] at (2,3.8) {\centerofmass};
      \draw[semithick,red!75!black,->,>=stealth]  (2,3.8) -- (2,2.8);
      \node[right] at (2,3.3) {\color{red!75!black}{$\mathrm{g}$}};
      
      \fill[
      left color=gray!50!black,
      middle color=gray!50,
      right color=gray!50!black,
      shading=axis,
      shading angle=35,
      opacity=0.25
      ] (3.75,3.75) -- (0,1.5) -- (1.3,5.2) [rotate=335] arc (180:340:1.49cm and 0.5cm); 
      \fill[
      shading=ball,
      ball color=gray!50,
      opacity=0.25
      ] (1.3,5.2) [rotate=335] arc (180:340:1.49cm and 0.5cm) arc (340:560:1.507cm); 
    \end{tikzpicture}
    \label{fig:htmb}
  }
  \caption{Motivating examples.}
  \label{fig:motivating_examples}
\end{figure}

The goal of this paper is to extend the method of controlled Lagrangians to a class of mechanical systems on semidirect product Lie groups with broken symmetry in order to find controls that stabilize their unstable equilibria.
Our motivating examples are a bottom-heavy underwater vehicle and a heavy top spinning on a movable base shown in \cref{fig:motivating_examples}.

These system, although seemingly quite different, have a few features in common:
\begin{enumerate}
\item Their configuration space is the semidirect product Lie group $\SE(3) \defeq \SO(3) \ltimes \R^{3}$.
\item One cannot decouple the dynamics into those in the rotational dynamics in  $\SO(3)$ and the translational dynamics in $\R^{3}$ as in the standard rigid body dynamics due to their interactions.
\item The gravity breaks their $\SE(3)$-symmetry the system would otherwise possess.
\end{enumerate}

Motivated by the first two features, we would like to consider mechanical systems whose configuration manifold is a semidirect product Lie group $\mathsf{S} \defeq \mathsf{G} \ltimes V$.
If the $\mathsf{S}$-symmetry \textit{were not} broken, the system \textit{would} possess $\mathsf{S}$-symmetry, and as a result, one would be able to write the equations of motion as the standard Euler--Poincar\'e equation (see, e.g., \cite[Chapter 13]{MaRa1999}) on the Lie algebra $\mathfrak{s} \defeq \mathfrak{g} \ltimes V$ of $\mathsf{S}$.
However, the broken symmetry mentioned in the last feature prevents one from performing such a symmetry reduction.

In order to remedy the broken $\mathsf{S}$-symmetry, one may introduce the so-called \textit{advected parameters} to the formulation to recover the $\mathsf{S}$-symmetry.
Assuming that the advected parameters live in the dual $X^{*}$ of a vector space $X$, the resulting Euler--Poincar\'e equations with advected parameters~\cite{HoMaRa1998a,CeHoMaRa1998} give differential equations on $\mathfrak{s} \times X^{*}$.

\subsection{Controlled Lagrangians}
The method of controlled Lagrangians was originally developed for those systems described by the Euler--Lagrange equations~\cite{OrPeNiSi1998,Ha1999,Ha2000,BlLeMa2000,BlChLeMa2001,Ch2008}, and was also applied to the standard Euler--Poincar\'e systems~\cite{BlKrMaSa1992,BlChLeMaWo2000,BlLeMa2001}.
We also note that there is the Hamiltonian version developed in \cite{Sc1986,OrScMaMa2001,FuSu2001,OrSpGoBl2002,BlOrVa2002} (see also \cite[\S12.3]{NiSc1990}); the two approaches are known to be equivalent for a certain class of systems~\cite{ChBlLeMaWo2002}.

We extend the method of controlled Lagrangians to the Euler--Poincar\'e equations \textit{with advected parameters} for those mechanical systems whose configuration manifold is a semidirect product Lie group $\mathsf{S} = \mathsf{G} \ltimes V$---with a particular interest in the case with $\mathsf{S} = \SE(3) = \SO(3) \ltimes \R^{3}$ motivated by the examples shown above.

The main advantages of the Euler--Poincar\'e equations with advected parameters are the following:
\begin{enumerate}[label=\arabic*.]
\item The equations of motion are defined on the vector space $\mathfrak{s} \times X^{*}$.
\item It does not directly involve parametrizations of the group $\mathsf{S}$ such as the Euler angles, which are known to cause difficulties in numerical computations~\cite{ZeLeBl2012,LeLeMc2017} because of coordinate singularities.
\item The kinetic energy is typically defined in terms of a quadratic form defined on the vector space.
\end{enumerate}

These features, particularly the last one, are particularly desirable for the kinetic shaping with the method of controlled Lagrangians because it boils down to considering a different quadratic form on the vector space.
In other words, the matching condition we seek here is less general than what is usually referred to as the matching condition (see, e.g., \citet{BlOrVa2002}) in which one obtains a PDE for the controlled Lagrangian.
We rather assume an ansatz for the controlled Lagrangian as in \cite{BlLeMa2000,BlChLeMa2001,BlLeMa2001} for the matching, and then perform a stability analysis to ensure the stabilization of the unstable equilibrium of interest in each specific case.

We note that \citet{ChMa2000,ChMa2004} achieved stabilization of the heavy top spinning \textit{on the ground} by using \textit{internal rotors} attached to the top.
This is also an example of the method of controlled Lagrangian applied to an Euler--Lagrange equations with advected parameters.
However, our second motivating example is different from theirs:
First, ours is the heavy top spinning on a \textit{movable base} as opposed to the ground; hence our configuration space is $\SE(3)$ as opposed to $\SO(3)$ in theirs.
Second, our control is applied as an \textit{external force} to the movable base as opposed to a torque applied to the top via internal rotors.

We also note that our result on the underwater vehicle is different from those of \citet{Le1997b,WoLe2002}.
Our present work mainly focuses on the kinetic shaping, whereas \citet{Le1997b} focuses on the potential shaping---the topic of our companion paper~\cite{CoOh2022a}.
\citet{WoLe2002} use torques by internal rotors, whereas our control involves external forces only.
Our setting is more amenable to those controls applied by, e.g., jets attached to the body.

\section{Semidirect Product Lie Groups}
\label{sec:semidirect_product}
We first give a brief summary of semidirect product Lie groups with a particular attention to $\SE(3)$.
This section overlaps with the companion paper \cite{CoOh2022a}, but is included for completeness as well as to set the notation.

\subsection{Semidirect Product Lie Groups and Lie Algebras}
Let $\mathsf{G}$ be a Lie group, $V$ be a vector space, and $\GL(V)$ be the set of all invertible linear transformations on $V$.
Let $\lambda\colon \mathsf{G} \to \GL(V)$ be a (left) representation of $\mathsf{G}$ on $V$, i.e., $\lambda(g_{1} g_{2}) = \lambda(g_{1}) \lambda(g_{2})$ for any $g_{1}, g_{2} \in \mathsf{G}$.
We then define the semidirect product Lie group $\mathsf{S} \defeq \mathsf{G} \ltimes V$ under the multiplication
\begin{align*}
  s_{1} \cdot s_{2}
  = (g_{1}, x_{1}) \cdot (g_{2}, x_{2})
  = (g_{1} g_{2}, \lambda(g_{1}) x_{2} + x_{1}).
\end{align*}

Let $\mathfrak{g}$ be the Lie algebra of $\mathsf{G}$.
Then the representation $\lambda$ induces the Lie algebra representation $\lambda'\colon \mathfrak{g} \to \gl(V)$ as follows:
\begin{equation*}
  \lambda'(\xi) v \defeq \dzero{t}{ \lambda(\exp(t\xi)) v } = \lambda^{i}_{\alpha k} \xi^{\alpha} v^{k} = \xi_{V}(v),
\end{equation*}
where $\xi_{V}$ is the infinitesimal generator on $V$ corresponding to $\xi \in \mathfrak{g}$.
Then we have the semidirect product Lie algebra $\mathfrak{s} \defeq \mathfrak{g} \ltimes V$ equipped with the commutator
\begin{equation}
  \label{eq:ad-sdp}
  \ad_{(\xi,v)}(\eta,w)
  \defeq [(\xi,v), (\eta,w)]
  = \parentheses{
    \ad_{\xi}\eta,\,
    \lambda'(\xi)w - \lambda'(\eta)v
  }.
\end{equation}

One may also fix $v \in V$ in $\lambda'(\xi) v$ to regard $\xi \mapsto \lambda'(\xi) v$ as a linear map $\lambda'_{v}\colon \mathfrak{g} \to V$, i.e.,
\begin{equation*}
  \lambda'_{v}(\xi) \defeq \lambda'(\xi) v = \parentheses{ \lambda^{i}_{\alpha k} v^{k} } \xi^{\alpha}.
\end{equation*}
Then its dual $(\lambda_{v}')^{*}$ defines the momentum map $\mathbf{J}\colon T^{*}V \cong V \times V^{*} \to \mathfrak{g}^{*}$ as follows:
\begin{align*}
  \ip{ \mathbf{J}(v,a) }{ \xi } = \ip{ (\lambda_{v}')^{*}a }{ \xi }
  = \ip{ a }{ \lambda_{v}'(\xi) }
  = a_{k} \lambda^{k}_{\alpha j} v^{j} \xi^{\alpha},
\end{align*}
which results in
\begin{equation}
  \label{eq:J}
  \mathbf{J}(v,a) = \lambda^{k}_{\alpha j} v^{j} a_{k}.
\end{equation}
This is nothing but the so-called diamond operator $\diamond\colon V \times V^{*} \to \mathfrak{g}^{*}$ (see \citet{HoMaRa1998a,CeHoMaRa1998} and \citet[\S7.5]{HoScSt2009}), i.e., $v \diamond a = \mathbf{J}(v,a)$.

Let us also find an expression for the dual $\lambda'(\xi)^{*}$ of $\lambda'(\xi)$:
\begin{equation*}
  \ip{ \lambda'(\xi)^{*} a }{ v } = \ip{ a }{ \lambda'(\xi) v }
  \iff
  \parentheses{ \lambda'(\xi)^{*} a }_{i} v^{i} = a_{k} \lambda^{k}_{\alpha i} \xi^{\alpha} v^{i},
\end{equation*}
which gives
\begin{equation}
  \label{eq:lambda'*}
  \lambda'(\xi)^{*} a = \lambda^{k}_{\alpha i} \xi^{\alpha} a_{k}.
\end{equation}

We may now write the coadjoint representation on the dual $\mathfrak{s}^{*}$ of $\mathsf{S}$ as follows:
\begin{equation}
  \label{eq:adstar}
  \ad_{(\xi,v)}^{*}(\mu,a)
  = \parentheses{
    \ad_{\xi}^{*}\mu - \mathbf{J}(v,a),\,
    \lambda'(\xi)^{*}a
  }
  = \parentheses{
    c^{\beta}_{\gamma\alpha}\xi^{\gamma}\mu_{\beta} - \lambda^{k}_{\alpha j} v^{j} a_{k},\,
    \lambda^{l}_{\beta i} \xi^{\beta} a_{l}
  }.
\end{equation}

\begin{example}[$\SE(3) = \SO(3) \ltimes \R^{3}$]
  \label{ex:SE3}
  Consider the representation $\lambda\colon \SO(3) \to \GL(\R^{3}) = \GL(3,\R)$ defined by the standard matrix-vector multiplication, i.e.,
  \begin{equation*}
    \lambda(R)\mathbf{x} = R \mathbf{x}.
  \end{equation*}
  Then we can define the special Euclidean group $\SE(3) \defeq \SO(3) \ltimes \R^{3}$ under the following group multiplication:
  \begin{equation*}
    (R_{1}, \mathbf{x}_{1}) \cdot (R_{2}, \mathbf{x}_{2})
    = (R_{1} R_{2}, R_{1} \mathbf{x}_{2} + \mathbf{x}_{1}).
  \end{equation*}
  Another way of looking at $\SE(3)$ is that it is a matrix group
  \begin{equation*}
    \SE(3) = \setdef{ (R,\mathbf{x}) \defeq \begin{bmatrix} R & \mathbf{x} \\ \mathbf{0}^{T} & 1 \end{bmatrix} }{ R \in \SO(3),\, \mathbf{x} \in \R^{3} }
  \end{equation*}
  under the standard matrix multiplication.
  One then sees that the left translation of $(\dot{R}, \dot{\mathbf{x}}) \in T_{(R,\mathbf{x})}\SE(3)$ to the Lie algebra $\se(3) \defeq T_{(I,\mathbf{0})}\SE(3)$ is
  \begin{equation}
    \label{eq:Omega_v}
    \mathsf{L}_{(R,\mathbf{x})^{-1}}(\dot{R}, \dot{\mathbf{x}})
    = (R^{-1} \dot{R},\, R^{-1} \dot{\mathbf{x}})
    \eqdef (\hat{\Omega}, \bv),
  \end{equation}
  where $\hat{\Omega} \in \so(3)$ is the body angular velocity and $\bv$ is the translational velocity in the body frame.
  Note that we may identify $\hat{\Omega} \in \so(3)$ with $\bOmega \in \R^{3}$ via the hat map defined as
  \begin{equation}
    \label{eq:hat_map}
    \hat{(\,\cdot\,)}\colon \R^{3} \to \so(3);
    \qquad
    \mathbf{a} \mapsto
    \hat{a} \defeq
    \begin{bmatrix}
      0 & -a_{3} & a_{2} \\
      a_{3} & 0 & -a_{1} \\
      -a_{2} & a_{1} & 0
    \end{bmatrix}.
  \end{equation}
  So we may use $(\bOmega, \bv) \in \R^{3} \times \R^{3}$ as coordinates for $\se(3)$.
  
  Then we have
  \begin{equation*}
    \lambda'(\hat{\Omega}) \bv = \lambda_{\bv}'(\hat{\Omega}) = \dzero{t}{ \exp( t\hat{\Omega} ) \bv }
    = \hat{\Omega} \bv
    = \bOmega \times \bv
    = {\eps^{i}}_{\alpha k} \Omega^{\alpha} v^{k}.
  \end{equation*}
  Therefore, we have $\lambda^{i}_{\alpha k} = {\eps^{i}}_{\alpha k}$, and so \eqref{eq:J} and \eqref{eq:lambda'*} give
  \begin{equation*}
    \mathbf{J}(\bv, \bP) = {\eps^{k}}_{\alpha j} v^{j} P_{k} = \bv \times \bP,
    \qquad
    \lambda'(\hat{\Omega})^{*} \bP = {{\eps}^{k}}_{\alpha i} \Omega^{\alpha} P_{k} = \bP \times \bOmega.
    \qquad
  \end{equation*}
  
  Note that, using the above identification of $\R^{3}$ with $\so(3)$, the structure constants satisfy $c^{\alpha}_{\beta\gamma} = {\eps^{\alpha}}_{\beta\gamma}$ as well.
  So, using \eqref{eq:adstar}, we may write the coadjoint representation as follows:
  \begin{align*}
    \ad_{(\bOmega, \bv)}^{*} (\bPi, \bP)
    = \parentheses{
    \bPi \times \bOmega + \bP \times \bv,\,
    \bP \times \bOmega
    }.
  \end{align*}
\end{example}

In \cref{sec:SE3_ltimes_R4}, we consider further semidirect products $\SE(3) \ltimes \R^{3}$ and $\SE(3) \ltimes \R^{4}$, which crop up in the formulations of our motivating examples.

\section{Euler--Poincar\'e Equation with Advected Parameters}
\label{sec:EPwithAdP}
\subsection{Recovering Broken Symmetry of Lagrangian}
Consider a mechanical system defined on a semidirect product Lie group $\mathsf{S} = \mathsf{G} \ltimes V$ with Lagrangian $L_{\Gamma_{0}}\colon T\mathsf{S} \to \R$ with parameters $\Gamma_{0} \in X^{*}$, where $X^{*}$ is the dual of a vector space $X$.
Specifically, we consider the Lagrangian of the following form:
\begin{equation*}
  L_{\Gamma_{0}}(s, \dot{s}) = \frac{1}{2} \metric{\dot{s}}{\dot{s}} - U_{\Gamma_{0}}(s),
\end{equation*}
where $\metric{\,\cdot\,}{\,\cdot\,}$ is a left-invariant metric on $T\mathsf{S}$, i.e., for any $s, s_{0} \in \mathsf{S}$ and any $\dot{s} \in T_{s}\mathsf{S}$,
\begin{equation*}
  \metric{T_{s}\mathsf{L}_{s_{0}}(\dot{s})}{T_{s}\mathsf{L}_{s_{0}}(\dot{s})} = \metric{\dot{s}}{\dot{s}},
\end{equation*}
where $\mathsf{L}$ stands for the left translation, i.e., $\mathsf{L}_{s_{0}}(s) = s_{0} s$ for any $s_{0}, s \in \mathsf{S}$, and $T\mathsf{L}$ is its tangent lift.
So the kinetic term is $\mathsf{S}$-invariant.

Suppose however that the potential is \textit{not} $\mathsf{S}$-invariant, i.e., there exist $s_{0}, s \in \mathsf{S}$ such that $U_{\Gamma_{0}}(s_{0} s) \neq U_{\Gamma_{0}}(s)$.
This breaks the $\mathsf{S}$-symmetry of the Lagrangian $L_{\Gamma_{0}}$.
We further suppose that we can fix this in the following way:
Define an extended potential $U\colon \mathsf{S} \times X^{*} \to \R$ so that $U(s, \Gamma_{0}) = U_{\Gamma_{0}}(s)$ for any $s \in \mathsf{S}$, and let $\kappa\colon \mathsf{S} \to \GL(X)$ be a representation of $\mathsf{S}$ on $X$, and $\kappa^{*}\colon \mathsf{S} \to \GL(X^{*})$ be the induced representation on the dual $X^{*}$.
We assume that we can find an appropriate $\kappa$ so that we can recover the $\mathsf{S}$-symmetry of the potential:
For any $s_{0}, s \in \mathsf{S}$ and any $\Gamma \in X^{*}$,
\begin{equation*}
  U\parentheses{ s_{0}s, \kappa(s_{0})^{*}\Gamma } = U(s, \Gamma).
\end{equation*}

Now let us define an extended Lagrangian $L\colon T\mathsf{S} \times X^{*} \to \R$ by setting
\begin{equation*}
  L(s, \dot{s}, \Gamma) \defeq \frac{1}{2} \metric{\dot{s}}{\dot{s}} - U(s, \Gamma),
\end{equation*}
and also define the action
\begin{align*}
  \Psi\colon &\mathsf{S} \times (T\mathsf{S} \times X^{*}) \to T\mathsf{S} \times X^{*};
  \\
  &(s_{0}, (s, \dot{s}, \Gamma)) \mapsto \Psi_{s_{0}}(s, \dot{s}, \Gamma) \defeq \parentheses{ s_{0} s, T_{s}\mathsf{L}_{s_{0}}(\dot{s}), \kappa^{*}(s_{0}) \Gamma }.
\end{align*}
Then we see that the extended Lagrangian now possesses the $\mathsf{S}$-symmetry, i.e., $L \circ \Psi_{s_{0}} = L$ for any $s_{0} \in \mathsf{S}$.

\subsection{Euler--Poincar\'e Equation with Advected Parameters}
Defining, with an abuse of notation, the reduced potential
\begin{equation*}
  U\colon X^{*} \to \R;
  \qquad
  U(\Gamma) \defeq U(e, \Gamma),
\end{equation*}
we may define the reduced extended Lagrangian $\ell\colon \mathfrak{s} \times X^{*} \to \R$ as
\begin{equation}
  \label{eq:ell}
  \ell(\xi, v, \Gamma) \defeq L(e, (\xi, v), \Gamma)
  = K(\xi, v) - U(\Gamma)
\end{equation}
with the kinetic energy term $K$ defined as
\begin{equation}
  \label{eq:K}
  K(\xi,v)
  \defeq \frac{1}{2} \metric{(\xi,v)}{(\xi,v)}
  = \frac{1}{2} \mathbb{G}_{\alpha\beta} \xi^{\alpha} \xi^{\beta} + \mathbb{G}_{\alpha j} \xi^{\alpha} v^{j} + \frac{1}{2} \mathbb{G}_{ij} v^{i} v^{j},
\end{equation}
where all these $\mathbb{G}$'s are constant matrices, and $\mathbb{G}_{\alpha\beta}$ and $\mathbb{G}_{ij}$ are assumed to be symmetric.
We also define $\mathbb{G}_{i\beta} \defeq \mathbb{G}_{\beta i}$ component-by-component so that $\mathbb{G}_{\alpha j} \xi^{\alpha} v^{j} = \mathbb{G}_{i\beta} v^{i} \xi^{\beta}$.
Then we obtain the Euler--Poincar\'e equation with advected parameters (see \cite{HoMaRa1998a,CeHoMaRa1998} and \cite[\S7.5]{HoScSt2009}):
\begin{equation*}
  \od{}{t}\parentheses{ \fd{\ell}{(\xi,v)} } = \ad_{(\xi,v)}^{*} \fd{\ell}{(\xi,v)} + \mathbf{K}\parentheses{ \fd{\ell}{\Gamma}, \Gamma },
  \qquad
  \od{\Gamma}{t} = \kappa'(\xi,v)^{*} \Gamma,
\end{equation*}
where we defined, for any smooth function $f\colon E \to \R$ on a real vector space $E$, its functional derivative $\tfd{f}{x} \in E^{*}$ at $x \in E$ such that, for any $\delta x \in E$, under the natural dual pairing $\ip{\,\cdot\,}{\,\cdot\,}\colon E^{*} \times E \to \R$,
\begin{equation*}
  \ip{ \fd{f}{x} }{ \delta x } = \dzero{t}{ f(x + t\,\delta x) }.
\end{equation*}
For example, if $E = \R^{n}$ with the dual pairing in terms of the dot product, $\tfd{f}{\mathbf{x}} = \tpd{f}{\mathbf{x}}$, i.e., the standard gradient.
Note also that $\mathbf{K} \colon X \times X^{*} \to \mathfrak{s}^{*} = \mathfrak{g}^{*} \times V^{*}$ is the momentum map associated with the above action $\kappa$ defined in a similar manner to $\mathbf{J}$:
\begin{equation*}
  \mathbf{K}(x, \Gamma)
  = \parentheses{ \mathbf{K}_{\mathfrak{g}^{*}}(x, \Gamma),\, \mathbf{K}_{V^{*}}(x, \Gamma) }
  \defeq (\kappa'_{x})^{*}\Gamma,
\end{equation*}
where we split the components of $\mathbf{K}$ into those in $\mathfrak{g}^{*}$ and $V^{*}$ as $\mathbf{K}_{\mathfrak{g}^{*}}$ and $\mathbf{K}_{V^{*}}$.
Then, using the formula~\eqref{eq:adstar} for the coadjoint action on $\mathfrak{s}^{*}$, we have
\begin{equation}
  \label{eq:EP}
  \begin{split}
    \od{}{t}\parentheses{ \fd{\ell}{\xi} } &= \ad_{\xi}^{*} \fd{\ell}{\xi} - \mathbf{J}\parentheses{v, \fd{\ell}{v}} + \mathbf{K}_{\mathfrak{g}^{*}}\parentheses{ \fd{\ell}{\Gamma}, \Gamma }, \\
    \od{}{t}\parentheses{ \fd{\ell}{v} } &= \lambda'(\xi)^{*} \fd{\ell}{v} + \mathbf{K}_{V^{*}}\parentheses{ \fd{\ell}{\Gamma}, \Gamma }, \\
    \od{\Gamma}{t} &= \kappa'(\xi,v)^{*} \Gamma.
  \end{split}
\end{equation}

\begin{example}[{Underwater vehicle~\cite{Le1997a,Le1997b,LeMa1997}; see also \cite{ChHaSmWi2007,SmChWiCa2009}}]
  \label{ex:uwv}
  Consider the underwater vehicle shown in \cref{fig:uwv-details}.
  The configuration space is $\mathsf{S} = \SE(3)$, i.e., rotations about the center of buoyancy and its translational positions.
  Let $\{ \mathbf{e}_{i} \}_{i=1}^{3}$ and $\{ \mathbf{E}_{i} \}_{i=1}^{3}$ be the orthonormal spatial/inertial and body frames, respectively.
  The orientation $R \in \SO(3)$ of the vehicle is defined so that $\mathbf{E}_{i} = R \mathbf{e}_{i}$ for $i = 1, 2, 3$.
  Note that our definitions of $\mathbf{e}_{3}$ and $\mathbf{E}_{3}$ are the opposite of those in \cite{Le1997a,Le1997b,LeMa1997}.
  Letting $\mathbf{x} \in \R^{3}$ be the position of the center of buoyancy in the spatial frame, we have an element $(R, \mathbf{x}) \in \SE(3)$ giving the orientation and the position of the vehicle.
  
  \begin{figure}[htbp]
    \centering
    \includegraphics[width=0.5\linewidth]{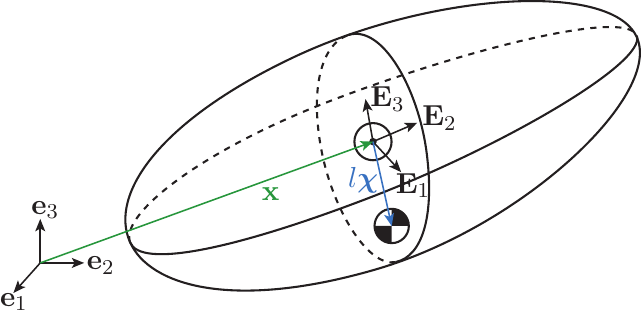}
    \caption{Underwater vehicle.}
    \label{fig:uwv-details}
  \end{figure}

  We assume that the vehicle is neutral buoyant and the shape of vehicle is ellipsoidal and also that the body frame is aligned with the principal axes of the body.
  Let $l \bchi$ be the position vector---$l$ being its length and $\bchi$ being the unit vector for the direction---of the center of mass measured from the center of buoyancy; see \cref{fig:uwv-details}.
  Then we have
  \begin{equation}
    \label{eq:mass_matrix-uwv}
    \mathbb{G}_{\alpha\beta} = \diag(I_{1}, I_{2}, I_{3}),
    \qquad
    \mathbb{G}_{\alpha j} = m l \hat{\chi},
    \qquad
    \mathbb{G}_{ij} = \diag(m_{1}, m_{2}, m_{3}).
  \end{equation}
  We note that $m_{1} \neq m_{2} \neq m_{3}$ in general, and so $\mathbb{G}_{ij}$ is not a constant multiple of the identity matrix; see \cite{Le1997a} for details.
   
  Due to the neutral buoyancy, the potential term is given as
  \begin{equation*}
    U_{\mathbf{e}_{3}}(R, \mathbf{x})
    = m \mathrm{g} l \mathbf{e}_{3} \cdot (R\bchi)
    = m \mathrm{g} l \bchi \cdot (R^{-1}\mathbf{e}_{3}).
  \end{equation*}
  Hence we define the extended potential $U\colon \SE(3) \times (\R^{3})^{*} \to \R$ by setting
  \begin{equation*}
    U((R, \mathbf{x}), \bGamma) \defeq m \mathrm{g} l \bchi \cdot (R^{-1}\bGamma)
  \end{equation*}
  so that $U((R, \mathbf{x}), \mathbf{e}_{3}) = U_{\mathbf{e}_{3}}(R, \mathbf{x})$.
 
  Using the representation~\eqref{eq:sigma-R3} of $\SE(3)$ on $\R^{3}$ from \cref{ssec:SE3-action_on_R3}, we have (see \eqref{eq:sigmas-R3}) 
  \begin{equation*}
    \kappa^{*}(R, \mathbf{x}) \bGamma = R \bGamma,
  \end{equation*}
  and so, for any $(R_{0}, \mathbf{x}_{0}), (R, \mathbf{x}) \in \SE(3)$ and any $\bGamma \in \R^{3}$,
  \begin{equation*}
    U\parentheses{ (R_{0}, \mathbf{x}_{0}) \cdot (R, \mathbf{x}), \kappa^{*}(R_{0}, \mathbf{x}_{0}) \bGamma }
    = U((R, \mathbf{x}), \bGamma).
  \end{equation*}
  Hence the reduced potential $U\colon (\R^{3})^{*} \to \R$ is
  \begin{equation*}
    U(\bGamma) \defeq U((I,\mathbf{0}),\bGamma) = m \mathrm{g} l \bchi \cdot \bGamma,
  \end{equation*}
  and the reduced Lagrangian $\ell\colon \se(3) \times (\R^{3})^{*} \to \R$ is
  \begin{equation*}
    \ell(\bOmega, \bv, \bGamma)
    = K(\bOmega,\bv) - U(\bGamma),
  \end{equation*}
  where $\bOmega$ and $\bv$ are defined in \eqref{eq:Omega_v}, and $K$ is the kinetic energy defined in \eqref{eq:K} using the mass matrix from \eqref{eq:mass_matrix-uwv}.
  Note that $\bGamma$ is the vertical upward direction (opposite of the direction of gravitational force) in the body frame.
  
  The representation $\kappa$ also gives (again see \eqref{eq:sigmas-R3})
  \begin{equation*}
    \kappa'(\bOmega,\bv)^{*} \bGamma = \bGamma \times \bOmega,
  \end{equation*}
  as well as the momentum map
  \begin{equation*}
    \mathbf{K}(\mathbf{y},\bGamma)
    = \parentheses{ \mathbf{K}_{\so(3)^{*}}(\mathbf{y},\bGamma), \mathbf{K}_{(\R^{3})^{*}}(\mathbf{y},\bGamma) }
    = (\kappa'_{\mathbf{y}})^{*}\bGamma
    = (\mathbf{y} \times \bGamma, \mathbf{0}).
  \end{equation*}
  As a result, the Euler--Poincar\'e equation~\eqref{eq:EP} with advected parameters gives
  \begin{equation}
    \label{eq:EP-uwv}
    \begin{split}
      \od{}{t} \parentheses{ \pd{\ell}{\bOmega} }
      &= \pd{\ell}{\bOmega} \times \bOmega 
      + \pd{\ell}{\bv} \times \bv
      + \pd{\ell}{\bGamma} \times \bGamma,
      \\
      \od{}{t} \parentheses{ \pd{\ell}{\bv} }
      &= \pd{\ell}{\bv} \times \bOmega,
      \\
      \dot{\bGamma}
      &= \bGamma \times \bOmega,
    \end{split}
  \end{equation}
  as in \cite{Le1997a,Le1997b,LeMa1997}.
\end{example}

\begin{example}[Heavy top on movable base]
  \label{ex:htmb}
  Consider the heavy top rotating on a movable base shown in \cref{fig:htmb-details}.
  \begin{figure}[hbtp]
    \centering
    \begin{tikzpicture}[scale=.75]
      \node [draw, cylinder, shape aspect=3, rotate=90, minimum height=1cm, minimum width=3cm] (c1)  at (0,1){};
      \draw[fill] (0,1.5) circle [radius=0.05];

      \fill[
      left color=gray!50!black,
      middle color=gray!50,
      right color=gray!50!black,
      shading=axis,
      opacity=0.25
      ] (3.75,3.75) -- (0,1.5) -- (1.3,5.2) [rotate=335] arc (180:340:1.49cm and 0.5cm); 
      \fill[
      shading=ball,
      ball color=gray!50,
      opacity=0.25
      ] (1.3,5.2) [rotate=335] arc (180:340:1.49cm and 0.5cm) arc (340:560:1.507cm); 

      \draw[semitransparent,fill] (2.5,3.5) circle [radius=0.025];
      \draw[semitransparent,->,>=stealth] (0,1.5) -- (2.5,3.5) node[right] {$\mathbf{X}$};
      \node at (1,5.5)    {$m$};
      \node at (2,0.5)    {$M$};

      \node[semitransparent] at (2,3.8) {\centerofmass};
      \draw[semithick,blue,->,>=stealth]  (0,1.5) -- (2,3.8);
      \node[left] at (1.9,3.7) {\color{blue}{$l\bchi$}};

      \draw[semithick,->,>=stealth]  (-5.5,1) -- (-5.5,1.5) node[above] {$\mathbf{e}_3$};
      \draw[semithick,->,>=stealth]  (-5.5,1) -- (-5.8,0.7) node[below] {$\mathbf{e}_1$};
      \draw[semithick,->,>=stealth]  (-5.5,1) -- (-5,1) node[below right] {$\mathbf{e}_2$};

      \draw[semithick,green!60!black][->,>=stealth]  (-5.5,1) -- (0,1.5);
      \node at (-3,1) {\color{green!60!black}{$\mathbf{x}$}};

      \draw[semithick,->,>=stealth]  (0,1.5) -- (-0.55,1.75) node[left] {$\mathbf{E}_1\!$};
      \draw[semithick,->,>=stealth]  (0,1.5) -- (0.29,1) node[right] {$\mathbf{E}_2$};
      \draw[semithick,->,>=stealth]  (0,1.5) -- (0.4,1.96) node[below right] {$\mathbf{E}_3$};

      \draw[very thick, ->,>=stealth, red!75!black] (-1.3,0.5) -- (-0.8,0.875);
      \node[below right, red!75!black] at (-1.3,0.575) {$\mathbf{u}$};

      \draw[->,>=stealth]  (-5.5,1) -- (2.5,3.5);
      \node[above] at (-1.5,2.25) {$q$};
    \end{tikzpicture}
    \caption{Heavy top on a movable base.}
    \label{fig:htmb-details}
  \end{figure}
  The configuration space is again $\SE(3)$:
  The orientation $R \in \SO(3)$ is defined in the same way as in \cref{ex:uwv} with respect to the body frame attached to the top at the junction point with the base, and is aligned with the principal axes; $\mathbf{x} \in \R^{3}$ is the position of the base, which is assumed to be a point mass $M$ for simplicity.
  
  Let $m$ be the mass of the heavy top, and $\bar{m} \defeq m+M$ the total mass of the system, $\mathbb{I} = \diag(I_{1},I_{2},I_{3})$ the inertia tensor of the top, $l$ the length of the line segment connecting the origin of the body frame (junction of body and base) to the center of mass of the top, $\bchi$ the unit vector pointing in that direction in the body frame, and $\mathrm{g}$ the gravitational constant---not to be confused with the italic $g$ used for an element of Lie group $\mathsf{G}$.
  
  Let $\mathcal{B} \subset \R^{3}$ be the domain occupied by the top in the body frame and $\rho_0\colon \mathcal{B} \to \R$ be the mass density of the top.
  Since the position $q(t)$ in the spatial frame of any point $\boldsymbol X\in \mathcal{B}$ at time $t$ is $q(t) = R(t) \boldsymbol X+ \mathbf{x}(t)$, the velocity of this point in the spatial frame is $\dot{q} = \dot{R} \boldsymbol X + \dot{\mathbf{x}}$.
  Therefore, we have the following Lagrangian:
  \begin{align*}
    \hspace{-0.25ex}
    L(R,\mathbf{x},\dot{R},\dot{\mathbf{x}})
    &= \int_{\mathcal{B}}
      \biggl(
      \frac{1}{2} \rho_0(\mathbf{X}) \| \dot{q} \|^2
    - \mathrm{g} \rho_0(\mathbf{X}) ( R\mathbf{X} + \mathbf{x} ) \cdot \mathbf{e}_3
      \biggr) d^3\mathbf{X}
      + \frac{1}{2} M  \| \dot{\mathbf{x}} \|^2
      - M\mathrm{g} \mathbf{x} \cdot \mathbf{e}_3
      \\
    &= \frac{1}{2}
      \biggl(
      \bar{m} \| \bv \|^2
      + \mathbb{I} \bOmega \cdot \bOmega
      + 2 \bv \cdot \bigl(\bOmega \times l \bchi m \bigr)
      \biggr)
      - m \mathrm{g} l \bchi \cdot (R^{-1} \mathbf{e}_{3})
      - \bar{m} \mathrm{g} \mathbf{x} \cdot \mathbf{e}_{3} \\
    &= K(\bOmega,\bv) - U_{e_{3}}(R, \mathbf{x})
  \end{align*}
  where the kinetic energy $K$ is defined in \eqref{eq:K} with
  \begin{equation}
    \label{eq:mass_matrix-htmb}
    \mathbb{G}_{\alpha\beta} = \diag(I_{1}, I_{2}, I_{3}),
    \qquad
    \mathbb{G}_{ij} = \bar{m} I,
    \qquad
    \mathbb{G}_{\alpha j} = m l \hat{\chi},
  \end{equation}
  and the potential term is defined as
  \begin{equation*}
    U_{e_{3}}(R, \mathbf{x})
    \defeq m \mathrm{g} l \bchi \cdot (R^{-1} \mathbf{e}_{3})
      + \bar{m} \mathrm{g} \mathbf{x} \cdot \mathbf{e}_{3}
   = \mathrm{g} \mathfrak{m} \cdot \parentheses{ s^{T} e_{3} }
  \end{equation*}
  with
  \begin{equation*}
    s = (R,\mathbf{x}) =
    \begin{bmatrix}
      R & \mathbf{x} \\
      \mathbf{0}^{T} & 1
    \end{bmatrix}
    \in \SE(3),
    \qquad
    \mathfrak{m} \defeq
    \begin{bmatrix}
      m l \bchi \\
      \bar{m}
    \end{bmatrix}
    \in \R^{4},
    \qquad
    e_{3} \defeq
    \begin{bmatrix}
      \mathbf{e}_{3} \\
      0
    \end{bmatrix}
    \in \R^{4}.
  \end{equation*}
  Notice that the potential $U_{e_{3}}$ depends not only on the orientation of the top but also on the height of the system, and hence is not $\SE(3)$-invariant.

  Let us define the extended potential $U\colon \SE(3) \times (\R^{4})^{*} \to \R$ by setting
  \begin{equation*}
    U((R, \mathbf{x}), \Gamma) \defeq \mathrm{g} \mathfrak{m} \cdot \parentheses{ s^{T} \Gamma },
  \end{equation*}
  so that $U((R, \mathbf{x}), e_{3}) = U_{e_{3}}(R, \mathbf{x})$.
  Using the representation $\kappa\colon \SE(3) \to \GL(\R^{4})$ defined in \eqref{eq:sigma-R4} in \cref{ssec:SE3-action_on_R4}, we have (see \eqref{eq:sigma*-R4})
  \begin{equation*}
    \kappa^{*}(s) \Gamma = (s^{T})^{-1} \Gamma.
  \end{equation*}
  As a result, we have, for any $s_{0}, s \in \SE(3)$,
  \begin{equation*}
    U\parentheses{ s_{0} s, \kappa^{*}(s_{0}) \Gamma }
    = U(s, \Gamma).
  \end{equation*}
  Let us write $\Gamma = (\bGamma, h) \in (\R^{4})^{*}$.
  Note that $\bGamma$ is the vertical upward direction (opposite of the direction of gravitational force) in the body frame, whereas $h$ is the height of the base in the inertial frame.
  Then we may define the reduced potential $U\colon (\R^{4})^{*} \to \R$ as
  \begin{equation}
    \label{eq:U-htmb}
    U(\bGamma, h) \defeq U(e, (\bGamma, h))
    = \mathrm{g} \mathfrak{m} \cdot (\bGamma, h)
    = m \mathrm{g} l \bchi \cdot \bGamma + \bar{m} \mathrm{g} h,
  \end{equation}
  and thus we have the reduced Lagrangian $\ell\colon \se(3) \times (\R^{4})^{*} \to \R$ as follows:
  \begin{equation*}
    \ell(\bOmega, \bv, (\bGamma, h))
    = K(\bOmega,\bv) - U(\bGamma, h).
  \end{equation*}
  
  Then, using the expressions in \eqref{eq:sigma'star-R4} and \eqref{eq:mathbfK-R4}, the Euler--Poincar\'e equation~\eqref{eq:EP} with advected parameters becomes:
  \begin{equation}
    \label{eq:EP-htmb}
    \begin{split}
      \od{}{t} \parentheses{ \pd{\ell}{\bOmega} }
      &= \pd{\ell}{\bOmega} \times \bOmega 
      + \pd{\ell}{\bv} \times \bv
      + \pd{\ell}{\bGamma} \times \bGamma,
      \\
      \od{}{t} \parentheses{ \pd{\ell}{\bv} }
      &= \pd{\ell}{\bv} \times \bOmega
      + \pd{\ell}{h}\bGamma,
      \\
      \dot{\bGamma} &= \bGamma \times \bOmega,
      \\
      \dot{h} &= \bGamma \cdot \bv.
    \end{split}
  \end{equation}
\end{example}

\begin{remark}
  \label{rem:potential_shaping}
  The above equations~\eqref{eq:EP-htmb} are very similar to \eqref{eq:EP-uwv} for the underwater vehicle.
  Indeed, one may apply the control force
  \begin{equation}
    \label{eq:u^p}
    \mathbf{u}^{\rm p} = -\pd{\ell}{h}\bGamma = \bar{m} \mathrm{g} \bGamma
  \end{equation}
  to the second equation of \eqref{eq:EP-htmb} to cancel the extra term, and as a result, may discard the height variable $h$ from the formulation to reduce the system to the same equation~\eqref{eq:EP-uwv} (with a slightly different kinetic energy metric).
  One can think of the above control as the potential shaping that cancels the second term on the right-hand side of \eqref{eq:U-htmb}; see the companion paper \cite{CoOh2022a} for details.
\end{remark}

\section{Controlled Lagrangian and Matching}
\subsection{Controlled Euler--Poincar\'e Equation with Advected Parameters}
Suppose that we would like to stabilize an unstable equilibrium of the system~\eqref{eq:EP} by applying an external (linear) force $u^{\rm k}$ (the roman superscript ``${\rm k}$'' denotes kinetic, not a coordinate index) to the system.
Practically speaking, the system is either pushed by some external means or controlled by jets attached to the body; the latter is more amenable to our formulation because our equations are written in the body frame.

Consider the controlled Euler--Poincar\'e equation with advected parameters:
\begin{equation}
  \label{eq:ControlledEP}
  \begin{split}
    \od{}{t}\parentheses{ \fd{\ell}{\xi} } &= \ad_{\xi}^{*} \fd{\ell}{\xi} - \mathbf{J}\parentheses{v, \fd{\ell}{v}} + \mathbf{K}_{\mathfrak{g}^{*}}\parentheses{ \fd{\ell}{\Gamma}, \Gamma }, \\
    \od{}{t}\parentheses{ \fd{\ell}{v} } &= \lambda'(\xi)^{*} \fd{\ell}{v} + \mathbf{K}_{V^{*}}\parentheses{ \fd{\ell}{\Gamma}, \Gamma } + u^{\rm k}, \\
    \od{\Gamma}{t} &= \kappa'(\xi,v)^{*} \Gamma.
  \end{split}
\end{equation}
We would like to match this control system with the Euler--Poincar\'e equation with advected parameters for a different reduced Lagrangian $\ellc\colon \mathfrak{s} \times X^{*} \to \R$:
\begin{equation}
  \label{eq:EP-ellc}
  \begin{split}
    \od{}{t}\parentheses{ \fd{\ellc}{\xi} } &= \ad_{\xi}^{*} \fd{\ellc}{\xi} - \mathbf{J}\parentheses{v, \fd{\ellc}{v}} + \mathbf{K}_{\mathfrak{g}^{*}}\parentheses{ \fd{\ellc}{\Gamma}, \Gamma }, \\
    \od{}{t}\parentheses{ \fd{\ellc}{v} } &= \lambda'(\xi)^{*} \fd{\ellc}{v} + \mathbf{K}_{V^{*}}\parentheses{ \fd{\ellc}{\Gamma}, \Gamma }, \\
    \od{\Gamma}{t} &= \kappa'(\xi,v)^{*} \Gamma.
  \end{split}
\end{equation}
In other words, we would like to find the controlled Lagrangian $\ellc$ such that \eqref{eq:EP-ellc} gives \eqref{eq:ControlledEP}.
Then we determine the control $u^{\rm k}$ such that \eqref{eq:ControlledEP} and \eqref{eq:EP-ellc} become equivalent.
As a result, the dynamics of the controlled system \eqref{eq:ControlledEP} is described by the ``free'' system \eqref{eq:EP-ellc} with the new Lagrangian $\ellc$.

\subsection{Controlled Lagrangian}
We would like to seek the controlled Lagrangian of the form
\begin{equation}
  \label{eq:ell-controlled}
  \ell_{\tau, \sigma, \rho}(\xi, v, \Gamma)
  \defeq K_{\tau, \sigma, \rho}(\xi, v) - U(\Gamma),
\end{equation}
where $K_{\tau, \sigma, \rho}$ is the modified kinetic energy whose expression we now seek in the following form as in \cite{BlLeMa2001}:
Using the kinetic energy $K$ and the metric tensor $\mathbb{G}$ from \eqref{eq:K} as well as constant matrices $\sigma$, $\rho$, and $\tau$ ($\sigma$ and $\rho$ are symmetric) to be determined below,
\begin{equation*}
  \begin{split}
    K_{\tau, \sigma, \rho}(\xi, v)
    &\defeq K\parentheses{ \xi^{\alpha}, v^{i} + \tau^{i}_{\alpha}\xi^{\alpha} }
    + \frac{1}{2} \sigma_{ij} \tau^{i}_{\alpha} \tau^{j}_{\beta} \xi^{\alpha} \xi^{\beta} \\
    &\quad + \frac{1}{2} (\rho_{ij} - \mathbb{G}_{ij}) \bigl( v^{i} + (\mathbb{G}^{ik}\mathbb{G}_{k\alpha} + \tau^{i}_{\alpha})\xi^{\alpha} \bigr) \bigl( v^{j} + (\mathbb{G}^{jk}\mathbb{G}_{k\beta} + \tau^{j}_{\beta})\xi^{\beta} \bigr) \\
    &= \frac{1}{2} ({\mathbb{G}}_{\alpha\beta} + \Delta_{\alpha\beta}) \xi^{\alpha} \xi^{\beta}
    + (\mathbb{G}_{i\beta} + \Delta_{i\beta}) v^{i} \xi^{\beta}
    + \frac{1}{2} \rho_{ij} v^{i} v^{j} \\
    &= K(\xi, v)
    + \frac{1}{2} \Delta_{\alpha\beta} \xi^{\alpha} \xi^{\beta}
    + \Delta_{i\beta} v^{i} \xi^{\beta} 
    + \frac{1}{2} \Delta_{ij} v^{i} v^{j}
  \end{split}
\end{equation*}
with
\begin{gather*}
  \Delta_{\alpha\beta} \defeq \bigl( \mathbb{G}_{i \beta} + \sigma_{ij} \tau^{j}_{\beta} \bigr) \tau^{i}_{\alpha}
  + \Delta_{i\beta} \bigl( \mathbb{G}^{ik}\mathbb{G}_{k\alpha} + \tau^{i}_{\alpha} \bigr),
  \qquad
  \Delta_{i\beta} \defeq \rho_{ij} \bigl( \mathbb{G}^{jk} \mathbb{G}_{k\beta} + \tau^{j}_{\beta} \bigr) - \mathbb{G}_{i\beta},
  \\
  \Delta_{ij} \defeq \rho_{ij} - \mathbb{G}_{ij},
\end{gather*}
where $\mathbb{G}^{ij}$ stands for the \textit{inverse} of the matrix $\mathbb{G}_{ij}$, and we use the same convention for other matrices too.

\subsection{Matching Condition}
Clearly $\tfd{\ell}{\Gamma} = \tfd{\ellc}{\Gamma}$, and so, in order to have a matching, it is sufficient to impose
\begin{equation}
  \label{eq:matching-I}
  \fd{\ellc}{\xi} = \fd{\ell}{\xi},
  \qquad
  \mathbf{J}\parentheses{v, \fd{\ell}{v} - \fd{\ellc}{v}} = 0.
\end{equation}
Then \eqref{eq:ControlledEP} and \eqref{eq:EP-ellc} match under the control $u^{\rm k}$ given as
\begin{equation*}
  u^{\rm k} = \od{}{t}\parentheses{ \fd{\ell}{v} - \fd{\ellc}{v} }
  - \lambda'(\xi)^{*} \parentheses{ \fd{\ell}{v} - \fd{\ellc}{v} }.
\end{equation*}

The first condition in \eqref{eq:matching-I} is equivalent to $\Delta_{\alpha\beta} \xi^{\beta} + \Delta_{\alpha j} v^{j} = 0$ for any $\xi \in \mathfrak{g}$ and any $v \in V$.
Hence this reduces to $\Delta_{\alpha\beta} = 0$ and $\Delta_{\alpha j} = 0$.
Then $\Delta_{i\beta} = 0$ as well, but then this gives
\begin{equation}
  \label{eq:matching_condition1}
  \tau^{i}_{\beta} = \big( \rho^{ij} - \mathbb{G}^{ij} \big)\mathbb{G}_{j\beta},
  \tag{MC1}
\end{equation}
whereas substituting $\Delta_{i\beta} = 0$ into $\Delta_{\alpha\beta} = 0$, we obtain
\begin{equation*}
  \bigl( \mathbb{G}_{i \beta} + \sigma_{ij} \tau^{j}_{\beta} \bigr) \tau^{i}_{\alpha} = 0.
\end{equation*}
We see that this is satisfied if $\mathbb{G}_{i \beta} + \sigma_{ij} \tau^{j}_{\beta} = 0$, but then this in turn is satisfied if
\begin{equation}
  \label{eq:matching_condition2}
  \sigma^{ij} = \mathbb{G}^{ij} - \rho ^{ij}.
  \tag{MC2}
\end{equation}

On the other hand, the second condition in \eqref{eq:matching-I} is written as, using \eqref{eq:J},
\begin{equation*}
  \lambda^{k}_{\alpha j} v^{j} (\Delta_{k\beta}\xi^{\beta} + \Delta_{kl} v^{l}) = 0.
\end{equation*}
Taking $\Delta_{i\beta} = 0$ and the expression for $\Delta_{kl}$ into account, we have
\begin{equation*}
  \lambda^{k}_{\alpha j} (\rho_{kl} - \mathbb{G}_{kl}) v^{j} v^{l} = 0.
\end{equation*}
Since this holds for any $v \in V$, it implies that $\lambda^{k}_{\alpha j} (\rho_{kl} - \mathbb{G}_{kl})$ is skew-symmetric with respect to the indices $(j,l)$, i.e.,
\begin{equation}
  \label{eq:matching_condition3}
  \lambda^{k}_{\alpha l} (\rho_{kj} - \mathbb{G}_{kj}) = -\lambda^{k}_{\alpha j} (\rho_{kl} - \mathbb{G}_{kl}).
  \tag{MC3}
\end{equation}

To summarize, we have the following:
\begin{theorem}
  \label{thm:matching}
  Under the matching conditions \eqref{eq:matching_condition1}---\eqref{eq:matching_condition3} and the control law
  \begin{equation*}
    u^{\rm k}_{i} = (\mathbb{G}_{ij} - \rho_{ij}) \dot{v}^{j} - \lambda^{j}_{\beta i} (\mathbb{G}_{jk} - \rho_{jk}) \xi^{\beta} v^{k},
  \end{equation*}
  the controlled Euler--Poincar\'e equations~\eqref{eq:ControlledEP} with advected parameters for the Lagrangian~\eqref{eq:ell} and the Euler--Poincar\'e equations~\eqref{eq:EP-ellc} with advected parameters for the controlled Lagrangian~\eqref{eq:ell-controlled} are equivalent.
\end{theorem}

\begin{remark}
  For implementation purposes, we may get rid of the acceleration $\dot{v}$ from the above feedback control law because we can rewrite \eqref{eq:EP-ellc} so that $(\dot{\xi}, \dot{v})$ is given in terms of functions of $(\xi, v, \Gamma)$; see \cref{ex:matching-SE3} below for an expression for the case with $\mathsf{S} = \SE(3)$.
\end{remark}

\begin{remark}
  Let us give an intuitive interpretation of the matching conditions.
  The conditions \eqref{eq:matching_condition1} and \eqref{eq:matching_condition2} imply that we ``reshape'' the kinetic energy by replacing the mass matrix $\mathbb{G}_{ij}$ by $\rho_{ij}$ only, i.e., no modifications of the other parts of the mass matrix.
  This intuitively makes sense because we are applying controls only to the ``translational'' part $V$.
  On the other hand, \eqref{eq:matching_condition3} imposes a restriction on the form of $\rho_{ij}$ to ensure that the interaction term between the ``rotational'' and ``translational'' parts ($\mathfrak{g}$ and $V$ respectively) matches with the original system.
  This also makes sense because their interactions are governed by the law of nature and should not be affected by the control.
\end{remark}

\begin{example}[$\mathsf{S} = \SE(3)$]
  \label{ex:matching-SE3}
  As seen in \cref{ex:SE3}, $\lambda^{i}_{\alpha k} = {\eps^{i}}_{\alpha k}$ in this case, and so the third matching condition~\eqref{eq:matching_condition3} becomes ${\eps^{k}}_{\alpha l} (\rho_{kj} - \mathbb{G}_{kj}) = -{\eps^{k}}_{\alpha j} (\rho_{kl} - \mathbb{G}_{kl})$.
  One may select $\rho$ so that $\rho_{ij} - \mathbb{G}_{ij}$ becomes a non-zero constant multiple of the identity matrix, i.e.,
  \begin{equation}
    \label{eq:rho-SE3}
    \rho_{ij} = \mathbb{G}_{ij} - \mathcal{K}\,\delta_{ij}
    \text{ for some }
    \mathcal{K} \in \R \backslash\{0\}.
  \end{equation}
  Then the above condition becomes ${\eps^{j}}_{\alpha l} = -{\eps^{l}}_{\alpha j}$, which is trivially satisfied.
  The feedback control then becomes
  \begin{equation}
    \label{eq:u^k-SE3}
    \mathbf{u}^{\rm k} = \mathcal{K} \parentheses{ \dot{\bv} + \bOmega \times \bv }.
  \end{equation}
  Note that, as mentioned in the above remark, one may replace the acceleration term $\dot{\bv}$ by a function of $(\bOmega, \bv, \bGamma)$ as follows:
  Using \eqref{eq:EP-ellc} along with the matching conditions, we have
  \begin{equation*}
    \begin{bmatrix}
      \dot{\boldsymbol{\Omega}} \\
      \dot{\mathbf{v}}
    \end{bmatrix}
    =
    \mathbb{G}^{-1}
    \begin{bmatrix}
      \DS \boldsymbol{\Pi} \times \boldsymbol\Omega 
      + \mathbf{P} \times \mathbf{v}
      - m \mathrm{g} l \boldsymbol{\chi} \times \boldsymbol\Gamma
      \medskip\\
      \DS \mathbf{P} \times \boldsymbol\Omega
    \end{bmatrix}
    \quad
    \text{with}
    \quad
    \boldsymbol{\Pi} \defeq \fd{\ellc}{\bOmega},
    \quad
    \mathbf{P} \defeq \fd{\ellc}{\bv}.
  \end{equation*}
\end{example}

\section{Stability Analysis}
\subsection{The Energy--Casimir Method}
\label{ssec:energy-Casimir}
We would like to establish the stability of equilibria of the systems from \cref{ex:uwv,ex:htmb} by constructing an appropriate Lyapunov function.
As mentioned in \cref{rem:potential_shaping}, the system from \cref{ex:htmb} after the ad-hoc potential shaping control~\eqref{eq:u^p} reduces to \eqref{eq:EP-uwv} in \cref{ex:uwv} with a slightly different Lagrangian.
Therefore, we may write down both systems under control force $\mathbf{u}^{\rm k}$ from \eqref{eq:u^k-SE3} via the kinematic shaping as
\begin{equation}
  \label{eq:ControlledEP-SE3}
  \begin{split}
    \od{}{t} \parentheses{ \pd{\ell}{\bOmega} }
    &= \pd{\ell}{\bOmega} \times \bOmega 
    + \pd{\ell}{\bv} \times \bv
    + \pd{\ell}{\bGamma} \times \bGamma,
    \\
    \od{}{t} \parentheses{ \pd{\ell}{\bv} }
    &= \pd{\ell}{\bv} \times \bOmega
    + \mathbf{u}^{\rm k},
    \\
    \dot{\bGamma}
    &= \bGamma \times \bOmega,
  \end{split}
\end{equation}
or equivalently
\begin{equation}
  \label{eq:ControlledEP2-SE3}
  \begin{split}
    \od{}{t} \parentheses{ \pd{\ellc}{\bOmega} }
    &= \pd{\ellc}{\bOmega} \times \bOmega 
    + \pd{\ellc}{\bv} \times \bv
    + \pd{\ellc}{\bGamma} \times \bGamma,
    \\
    \od{}{t} \parentheses{ \pd{\ellc}{\bv} }
    &= \pd{\ellc}{\bv} \times \bOmega,
    \\
    \dot{\bGamma}
    &= \bGamma \times \bOmega,
  \end{split}
\end{equation}
with the controlled Lagrangian
\begin{equation}
  \label{eq:ellc-SE3}
  \ellc(\bOmega, \bv, \bGamma)
  \defeq K_{\tau, \sigma, \rho}(\bOmega, \bv) - U(\bGamma).
\end{equation}

The main advantage of the method of controlled Lagrangians is that, thanks to the matching, the controlled system possesses invariants (conserved quantities) such as the energy and Casimirs, and is amenable to the \textit{energy--Casimir method} (see, e.g., \cite[\S1.7]{MaRa1999}).
Its main idea is to use such invariants to construct an invariant $\mathcal{E}$ that works as a control Lyapunov function to establish the stability of the equilibrium.

More specifically, the energy--Casimir method also prescribes a method to find such a Lyapunov function $\mathcal{E}$ using the energy of the system as well as Casimirs (or some other invariants of the system) as follows:
It is straightforward to show that the energy
\begin{equation*}
  E_{\tau,\sigma,\rho}(\bOmega, \bv, \bGamma)
  \defeq K_{\tau, \sigma, \rho}(\bOmega, \bv) + U(\bGamma)
\end{equation*}
associated with the controlled Lagrangian~\eqref{eq:ellc-SE3} is an invariant of the system~\eqref{eq:ControlledEP2-SE3}.
Also, as mentioned in \cref{ssec:LPB}, the system~\eqref{eq:ControlledEP2-SE3} has three Casimir functions (see \eqref{eq:Casimirs-LP}) or in the Lagrangian variables,
\begin{equation}
  \label{eq:Casimirs}
  C_{1} = \norm{ \pd{\ell_{\tau, \sigma, \rho}}{\bv} }^{2},
  \qquad
  C_{2} = \pd{\ell_{\tau, \sigma, \rho}}{\bv} \cdot \bGamma,
  \qquad
  C_{3} = \norm{ \bGamma }^{2}.
\end{equation}
This implies that, for any smooth function $\Phi\colon \R^{3} \to \R$, the function
\begin{equation*}
  \mathcal{E} \defeq E_{\tau,\sigma,\rho} + \Phi(C_{1}, C_{2}, C_{3})  
\end{equation*}
is also an invariant of the system~\eqref{eq:ControlledEP2-SE3} as well.
Note that the actual form of $\mathcal{E}$ varies depending on whether the system has other invariants, as we shall see below.

Now, one determines $\Phi$ so that $\mathcal{E}$ provides a control Lyapunov function.
Specifically, let $\zeta_{\rm e}$ be an equilibrium of the uncontrolled system~\eqref{eq:EP-uwv}, and proceed as follows:
\begin{enumerate}[label=\arabic*.]
\item Find the conditions under which the first variation (the gradient) $D\mathcal{E}$ vanishes at the equilibrium $\zeta_{\rm e}$.
\item Calculate the second variation (the Hessian) $D^{2}\mathcal{E}$ at $\zeta_{\rm e}$.
\item Find the conditions under which the Hessian $D^{2}\mathcal{E}(\zeta_{\rm e})$ is definite.
\end{enumerate}
As a result, there exists an open neighborhood $U$ of $\zeta_{\rm e}$ such that $\mathcal{E}(\zeta) > \mathcal{E}(\zeta_{\rm e})$ (or $\mathcal{E}(\zeta) < \mathcal{E}(\zeta_{\rm e})$) for any $\zeta \in U \backslash \{ \zeta_{\rm e} \}$.
Note also that $\zeta_{\rm e}$ is an equilibrium of the controlled system \eqref{eq:ControlledEP2-SE3} as well because $\mathcal{E}$ is an invariant of \eqref{eq:ControlledEP2-SE3} and $\mathcal{E}(\zeta_{\rm e})$ is a strict local extremum.

As a result, $\mathcal{E}$ gives a control Lyapunov function, and hence Lyapunov's Stability Theorem (see, e.g., \citet[Theorem~4.1]{Kh2002} and \citet[Theorem~5.2]{LoRy2014}) implies that the equilibrium $\zeta_{\rm e}$ is stable.

\subsection{Heavy top on movable base}
Consider \cref{ex:htmb} (see also \cref{fig:htmb}) with the Lagrange top, i.e., the inertia tensor $\mathbb{I} = \diag(I_{1},I_{2},I_{3})$ satisfies $I_1 = I_2 \neq I_3$, and its center of mass lies on the axis of symmetry with respect to the body frame, that is, $\bchi= (0, 0, 1)$.
We would like to show that the top spinning upright on the stationary base can be stabilized by the above control.
Note that, combining $\mathbb{G}_{ij} = \bar{m}\,\delta_{ij}$ from \eqref{eq:mass_matrix-htmb} and $\rho_{ij} = \mathbb{G}_{ij} - \mathcal{K}\,\delta_{ij}$ from \eqref{eq:rho-SE3}, we may set $\rho_{ij} = \varrho\,\delta_{ij}$ with $\varrho \defeq \bar{m} - \mathcal{K} \in \R$.

This system has two additional invariants besides the energy and the Casimirs:
The first one is the well-known invariant $\Omega_{3}$ for the Lagrange top, and the second and less obvious one is the energy-like invariant:
\begin{equation*}
  E^{0}(\bOmega,\bGamma)
  \defeq \frac{1}{2}\parentheses{ I_{1}(\Omega_{1}^{2} + \Omega_{2}^{2}) + I_{3}\Omega_{3}^{2} }
  + m \frac{I_{1} \varrho}{I_{1} \varrho - m^{2} l^{2}} \mathrm{g} l \Gamma_{3}.
\end{equation*}
This implies that, for any constant $c \in \R$ and any smooth functions $\Phi\colon \R^{3} \to \R$ and $\phi\colon \R \to \R$,
\begin{equation}
  \label{eq:mathcalE-htmb}  
  \mathcal{E} \defeq
  E_{\tau,\sigma,\rho}
  + c\,E^{0}
  + \Phi(C_{1}, C_{2}, C_{3})
  + \phi(\Omega_{3})
\end{equation}
is also an invariant of the system as well.

The equilibrium corresponding to the top spinning upright on the stationary base is
\begin{equation}
  \label{eq:equilibrium-htmb}
  (\bOmega_{\text{e}}, \bv_{\text{e}}, \bGamma_{\text{e}}) = ( \Omega_{0} \mathbf{E}_{3}, \mathbf{0}, \mathbf{E}_{3} ).  
\end{equation}
Note that the upright spinning Lagrange top with $| \Omega_{0} | > 2\sqrt{m\mathrm{g} l I_{1}}/I_{3}$ is known to be stable~\cite[Theorem~15.10.1]{MaRa1999}.
Therefore we assume that $| \Omega_{0} | < 2\sqrt{m\mathrm{g} l I_{1}}/I_{3}$ here, and show that the equilibrium is stabilized regardless of the value of $\Omega_{0}$.

An interesting observation is that the above energy-like invariant $E^{0}$ is the energy of the Lagrange top \textit{without} the movable base---the only difference is that the gravitational constant $\mathrm{g}$ is modified to be $\frac{I_{1} \varrho}{I_{1} \varrho - m^{2} l^{2}} \mathrm{g}$.
This observation suggests the following: If we pick $\varrho \in (0, m^{2}l^{2}/I_{1})$, then the modified gravitational constant $\frac{I_{1} \varrho}{I_{1} \varrho - m^{2} l^{2}} \mathrm{g}$ becomes negative, and hence effectively turning the upright position of the top into the vertical downward one for the controlled system.
As a result, the upright position of the controlled system becomes stable.
Let us justify this intuitive argument using the energy--Casimir method.
\begin{proposition}[Stabilization of heavy top on a movable base]
  \label{prop:stabilization-htmb}
  The unstable equilibria~\eqref{eq:equilibrium-htmb} with $| \Omega_{0} | < 2\sqrt{m\mathrm{g} l I_{1}}/I_{3}$ of the heavy-top-on-movable-base system in \cref{ex:htmb} are stabilized by applying to the second equation of \eqref{eq:EP-htmb} the control $\mathbf{u} = \mathbf{u}^{\rm p} + \mathbf{u}^{\rm k}$, where $\mathbf{u}^{\rm p}$ is defined in \eqref{eq:u^p} and $\mathbf{u}^{\rm k}$ is from \eqref{eq:u^k-SE3} with $\mathcal{K} = \bar{m} - \varrho$ for any $\varrho \in (0,m^{2}l^{2}/I_{1})$.
\end{proposition}

\begin{proof}
  Note first that, as mentioned above, we have $\rho_{ij} = \varrho\, \delta_{ij}$ with $\varrho \defeq \bar{m} - \mathcal{K}$ here.
  
  Let us use $(\,\cdot\,)|_{\rm e}$ to indicate that a function is evaluated at the equilibrium.
  The first variation condition $D\mathcal{E}|_{\rm e} = 0$ is satisfied if
  \begin{equation}
    \label{eq:first_variation-htmb}
    D_{2}\Phi|_{\rm e} = 0,
    \qquad
    D_{3}\Phi|_{\rm e} = \frac{m^{2}l^{2} - (1+c)I_{1}\varrho}{2(I_{1} \varrho - m^{2} l^{2})} m g l,
    \qquad
    \phi'(\Omega_{0}) = -(1 + c) I_3 \Omega_{0},
  \end{equation}
  where $D_{i}$ stands for the derivative with respect to the $i$-th variable.
 
  By evaluating the leading principal minors of the Hessian $D^{2}\mathcal{E}|_{\rm e}$, we also find that the following conditions---in addition to \eqref{eq:first_variation-htmb}---are sufficient for its positive-definiteness:
  \begin{equation}
    \label{eq:posdef-htmb}
    \begin{array}{c}
      \DS D_{1}\Phi|_{\rm e} = D^{2}_{22}\Phi|_{\rm e} = D^{2}_{33}\Phi|_{\rm e} = D^{2}_{23}\Phi|_{\rm e} = \phi''(\Omega_{0}) = 0,
      \medskip\\
      c > 0,
      \qquad
      \DS \frac{m^{2}l^{2}}{(1 + c)I_{1}} < \varrho < \frac{m^{2}l^{2}}{I_{1}}.
    \end{array}
  \end{equation}
  Therefore, we may take, for example,
  \begin{equation}
    \label{eq:Phi-htmb}
    \Phi(C_{1}, C_{2}, C_{3}) = \frac{m^{2}l^{2} - (1+c)I_{1}\varrho}{2(I_{1} \varrho - m^{2} l^{2})} m g l (C_{3} - C_{3}|_{\rm e}),
    \quad
    \phi(\Omega_{3}) = - (1 + c) I_{3} \Omega_{0} (\Omega_{3} - \Omega_{0}).
  \end{equation}
   
  However, since we may take $c > 0$ arbitrarily large, we can achieve stability for any $\varrho \in (0,m^{2}l^{2}/I_{1})$.
\end{proof}

\begin{figure}[htbp]
  \centering
  \subfloat[Body angular velocity $(\Omega_{1},\Omega_{2})$]{
    \includegraphics[width=.3\linewidth]{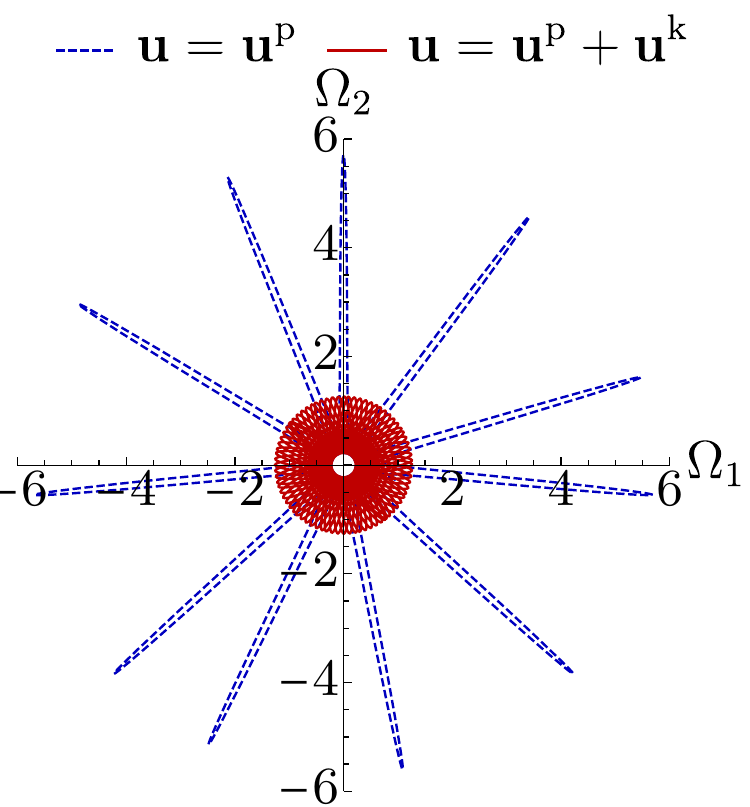}
  }
  \quad
  \subfloat[Base velocity $\bv$ in body frame]{
    \includegraphics[width=.29\linewidth]{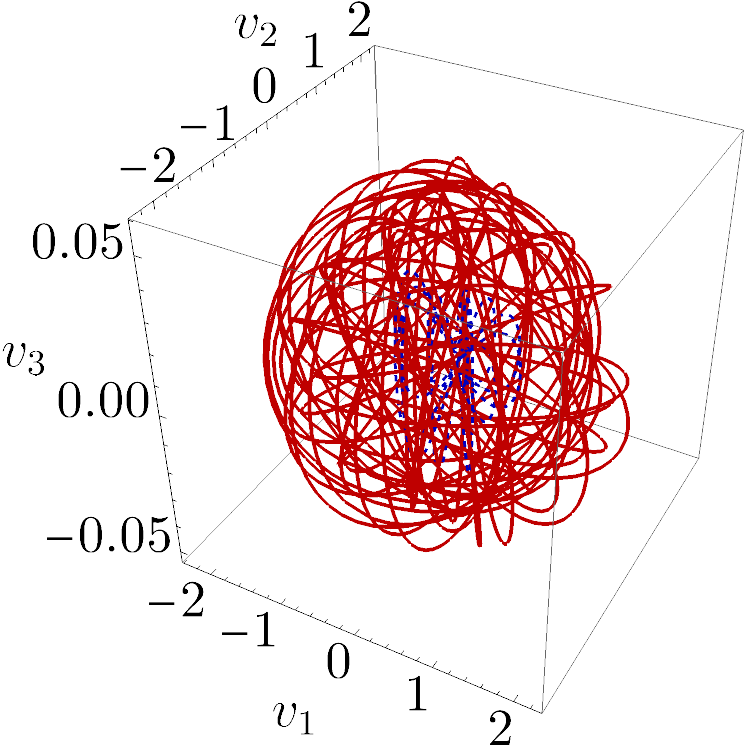}
  }
  \quad
  \subfloat[Vertical upward direction $\boldsymbol{\Gamma}$ seen from body frame]{
    \includegraphics[width=.3\linewidth]{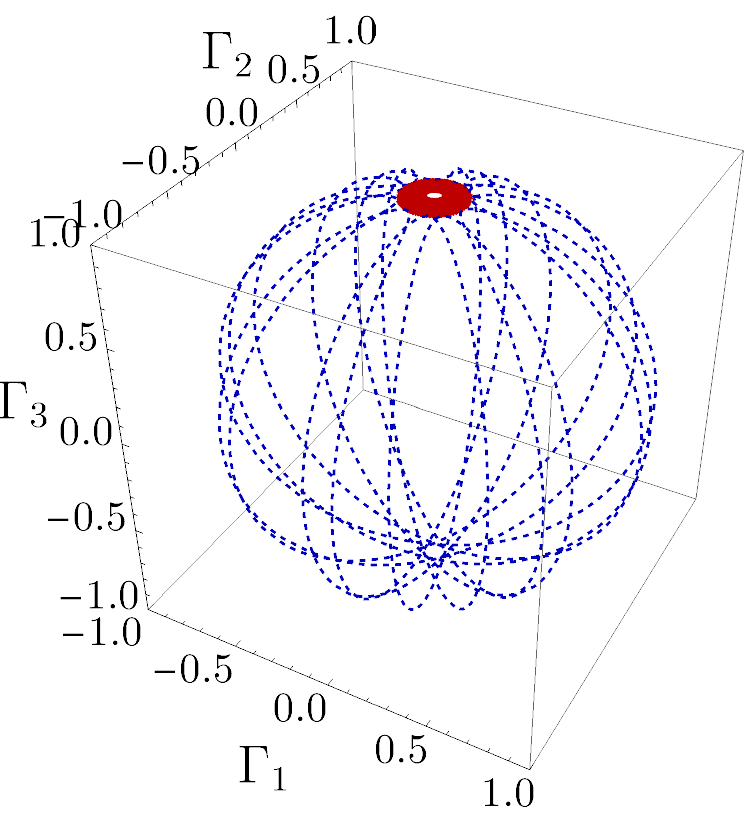}
  }
  \caption{Simulation results for the heavy top on a movable base with $M=0.44\,\mathrm{[kg]}$, $m=0.7\,\mathrm{[kg]}$, $I_{1} = I_{2} = 0.2\,\mathrm{kg \cdot m^{2}}$, $I_{3} = 0.24\,\mathrm{kg \cdot m^{2}}$, $l = 0.215\,\mathrm{[m]}$, $\mathrm{g} = 9.8\,\mathrm{[m/s^{2}]}$, and $\varrho = 0.9m^{2}l^{2}/I_{1}$ with initial condition $\bOmega(0) = (0.1, 0.2, 0.1)$, $\bv(0) = \mathbf{0}$, and $\boldsymbol{\Gamma}(0) = (\cos\theta_{0} \sin\varphi_{0}, \sin\theta_{0} \sin\varphi_{0}, \cos\varphi_{0})$ with $\theta_{0} = \pi/3$ and $\varphi_{0} = \pi/20$.
    The solutions are shown for the time interval $0 \le t \le 30$.
    The blue dashed line is for the system with only the control \eqref{eq:u^p} coming from potential shaping, whereas the red solid line is for the system with both the potential and kinetic shaping controls \eqref{eq:u^p} and \eqref{eq:u^k-SE3}.
    Note that the uncontrolled system with $\mathbf{u} = \mathbf{0}$ involves a free fall and does not provide a good comparison to illustrate the effect of stabilizing control $\mathbf{u}^{\rm k}$.}
  \label{fig:simulation-HT}
\end{figure}

\Cref{fig:simulation-HT} shows the results of simulations demonstrating the stabilizing control by the kinetic shaping; see the caption for the parameters and initial condition.
One can see that the equilibrium~\eqref{eq:equilibrium-htmb} is unstable without control $\mathbf{u}^{\rm k}$, but is stabilized after the control is applied to the system.

\subsection{Underwater Vehicle}
Let us now consider the underwater vehicle from \cref{ex:uwv}.
We assume, in addition to those assumptions mentioned in \cref{ex:uwv}, that the center of mass is aligned with the third principal axis $\mathbf{E}_{3}$ and below the center of buoyancy, i.e., $\bchi = (0,0,-1)$, and so it is bottom-heavy.

The equilibrium of our interest is the steady translational motion along $\mathbf{E}_{2}$, i.e.,
\begin{equation}
  \label{eq:equilibrium-uwv}
  (\bOmega_{\rm e}, \mathbf{v}_{\rm e}, \bGamma_{\rm e})
  = (\mathbf{0}, v_{0}\mathbf{E}_{2}, \mathbf{E}_{3})
\end{equation}
with $v_{0} \in \R\backslash\{0\}$.
According to \citet[Theorem~2]{Le1997a}, this is an unstable equilibrium of \eqref{eq:EP-uwv} if the vehicle is bottom-heavy and $m_{2} < m_{1}$, which is the case if the semi-major axis of the ellipsoidal hull along $\mathbf{E}_{2}$ is longer than that along $\mathbf{E}_{1}$ as depicted in \cref{fig:uwv-details}; see \cite[Appendix~B]{Le1997a} for details.

We can show that our control~\eqref{eq:u^k-SE3} stabilizes this equilibrium too:
\begin{proposition}[Stabilization of underwater vehicle]
  \label{prop:stabilization-uwv}
  The unstable equilibria~\eqref{eq:equilibrium-uwv} with $v_{0} \in \R\backslash\{0\}$ of the underwater vehicle system~\eqref{eq:EP-uwv} in \cref{ex:uwv} are stabilized by applying to the second equation of \eqref{eq:EP-uwv} the control $\mathbf{u} = \mathbf{u}^{\rm k}$ (as in \eqref{eq:ControlledEP-SE3}), where $\mathbf{u}^{\rm k}$ is from \eqref{eq:u^k-SE3} with any $\mathcal{K}$ satisfying
  \begin{equation}
    \label{eq:mathcalK-stability}
    m_{2} < \mathcal{K} < \min\braces{m_{3}, m_{1} - \frac{m l^{2}}{I_{2}}m }.
  \end{equation}
\end{proposition}

\begin{proof}
  Let us seek the control Lyapunov function of the form
  \begin{equation}
    \label{eq:mathcalE-uwv}
    \mathcal{E}(\bOmega, \bv, \bGamma) \defeq E_{\tau, \sigma, \rho}(\bOmega, \bv, \bGamma) + \Phi(C_{1}, C_{2}, C_{3}),
  \end{equation}
  because this system does not seem to have any additional invariants besides the energy and the Casimirs.
  
  One can show that $D\mathcal{E}|_{\rm e} = 0$ if
  \begin{equation}
    \label{eq:first_variation-uwv}
    D_{1}\Phi|_{\rm e} = \frac{1}{2(\mathcal{K} - m_{2})},
    \qquad
    D_{2}\Phi|_{\rm e} = 0,
    \qquad
    D_{3}\Phi|_{\rm e} = \frac{m g l}{2}.
  \end{equation}
  On the other hand, by evaluating the leading principal minors of the Hessian $D^{2}\mathcal{E}|_{\rm e}$, one can show that it is positive-definite if, in addition to \eqref{eq:first_variation-uwv}, all the components of the Hessian of $\Phi$ vanish except
  \begin{equation*}
    D^{2}_{11}\Phi|_{\rm e} = \frac{1}{(\mathcal{K} - m_{2})^{3} v_{0}^{2}},
  \end{equation*}
  and also the parameter $\mathcal{K}$ satisfies \eqref{eq:mathcalK-stability}.
  
  This implies that one may take, e.g.,
  \begin{equation}
    \label{eq:Phi-uwv}
    \Phi(C_{1}, C_{2}, C_{3}) = \frac{1}{2} \parentheses{ \frac{( C_{1} - C_{1}|_{\rm e} )^{2}}{(\mathcal{K} - m_{2})^{3} v_{0}^{2}}  + \frac{C_{1} - C_{1}|_{\rm e}}{\mathcal{K} - m_{2}} + m g l (C_{3} - C_{3}|_{\rm e}) }
  \end{equation}
  to satisfy the above conditions.
\end{proof}

\begin{remark}
  There must exist $\mathcal{K}$ satisfying \eqref{eq:mathcalK-stability} for those underwater vehicles of interest here.
  In fact, one can show that $m_{2} < m_{1}$ and $m_{2} < m_{3}$ if the semi-major axis of the ellipsoidal hull along $\mathbf{E}_{2}$ is longer than those along $\mathbf{E}_{1}$ and $\mathbf{E}_{3}$ as depicted in \cref{fig:uwv-details}; see \cite[Appendix~B]{Le1997a}.
  We would also have $m_{i} > m$ for any $i = 1, 2, 3$ (again see \cite[Appendix~B]{Le1997a}) and $m l^{2}/I_{2} \ll 1$, and so $m_{2} < m_{1} - \frac{m l^{2}}{I_{2}}m$.
\end{remark}

As a numerical example, consider an underwater vehicle whose hull is an ellipsoidal shell with the outer semi-major axes $(a_{1}, a_{2}, a_{3}) = (5, 10, 4)\,\mathrm{[m]}$ and the inner semi-major axes $(a_{1} - h, a_{2} -h, a_{3} - h)$ with $h \simeq 0.1666\,\mathrm{[m]}$ made of steel with density $\num{8000}\,\mathrm{[kg/m^{3}]}$.
For simplicity, we assume that all extra weight is concentrated at the point 1 meter below the center of the ellipsoids as a point mass with $40\%$ of the weight of the shell; hence the center of mass is at $l \bchi$ with $l = 2/7\,\mathrm{[m]}$ and $\bchi = (0, 0, -1) = -\mathbf{E}_{3}$.
Then the total mass of the vehicle is $m = \num{835245}\,\mathrm{[kg]}$, and it is neutrally buoyant assuming that the mass density of the water is $997\,\mathrm{[kg/m^{3}]}$---the ``thickness'' $h$ of the hull is determined that way.
Using formulas from \cite[Appendix~B]{Le1997a}, one obtains $(m_{1}, m_{2}, m_{3}) \simeq (1.330, 0.9860, 1.592) \times 10^{6}\,\mathrm{[kg]}$ and $(I_{1}, I_{2}, I_{3}) \simeq (2.787, 0.9020, 2.527)  \times 10^{7}\,\mathrm{[kg \cdot m^{2}]}$.
We set $\mathcal{K} \simeq 1.239 \times 10^{6}$ so that \eqref{eq:mathcalK-stability} is satisfied.

We select an initial condition with a small perturbation to the equilibrium~\eqref{eq:equilibrium-uwv} with $v_{0} = 30\,\mathrm{[m/s]}$ as follows:
\begin{gather*}
  \bOmega(0) = (0.5,\, 0.25,\, 0.5),
  \qquad
  \bv(0) = (1.5,\, 30,\, 1.5),
  \\
  \bGamma(0) = (\cos\theta_{0} \sin\varphi_{0}, \sin\theta_{0} \sin\varphi_{0}, \cos\varphi_{0})
\end{gather*}
with $\theta_{0} = \pi/3$ and $\varphi_{0} = \pi/40$.

\begin{figure}[htbp]
  \centering
  \subfloat[Body angular velocity $\bOmega$]{
    \includegraphics[width=.31\linewidth]{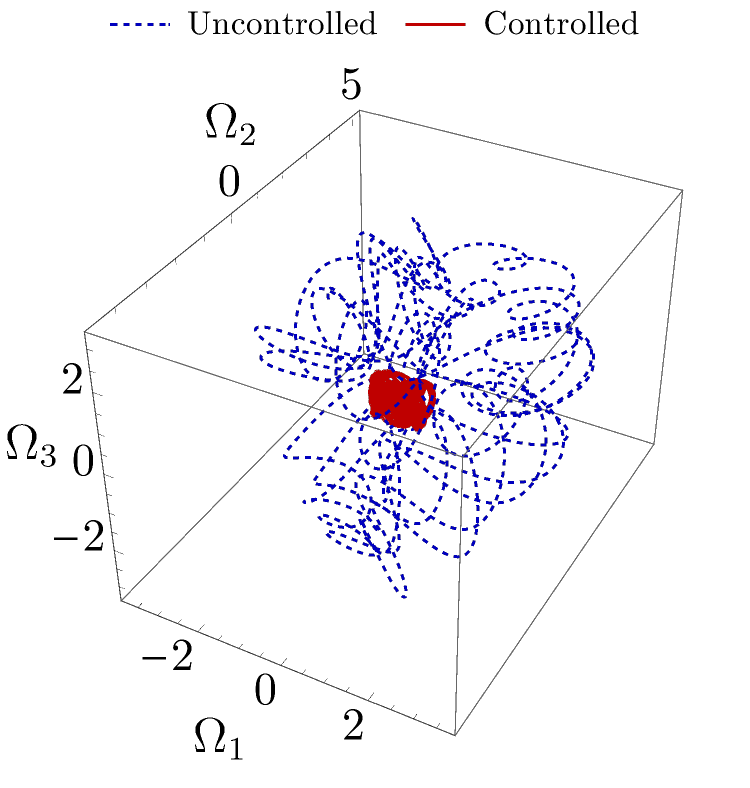}
  }
  \
  \subfloat[Base velocity $\mathbf{v}$ in body frame]{
    \includegraphics[width=.31\linewidth]{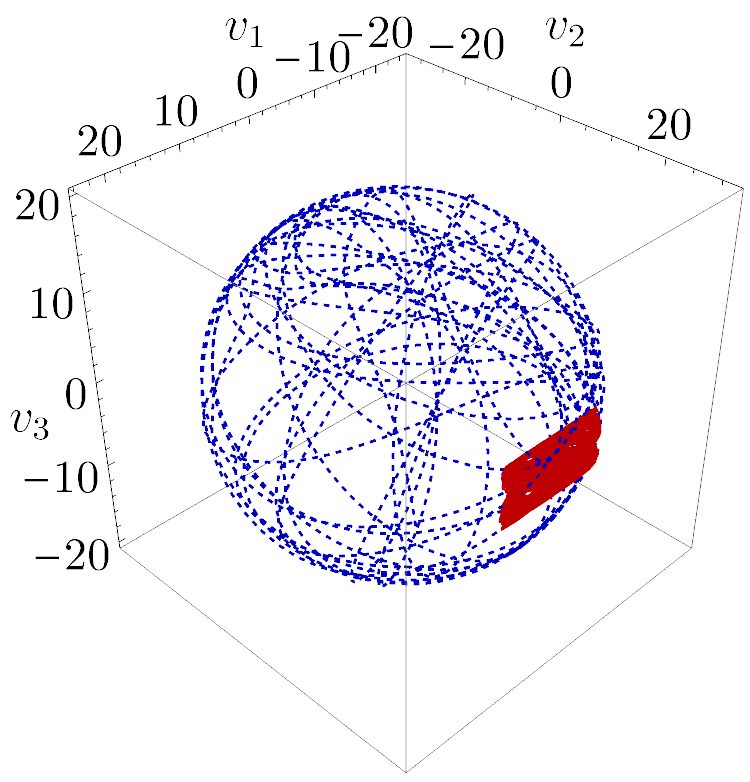}
  }
  \
  \subfloat[Vertical upward direction $\boldsymbol{\Gamma}$ seen from body frame]{
    \includegraphics[width=.31\linewidth]{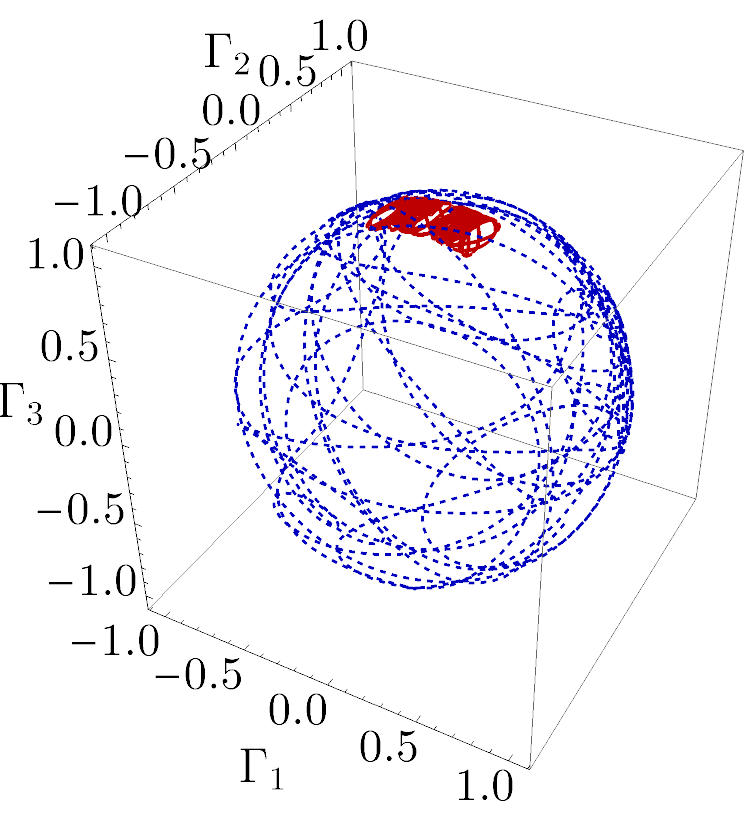}
  }
  \caption{Simulation results comparing the uncontrolled underwater vehicle~\eqref{eq:EP-uwv} (blue dashed) and controlled underwater vehicle~\eqref{eq:ControlledEP-SE3} (red solid) for the time interval $0 \le t \le 50$.}
  \label{fig:simulation-UWV}
\end{figure}

\Cref{fig:simulation-UWV} shows the trajectories of $\bOmega$, $\bv$, and $\bGamma$ for the uncontrolled and controlled systems.
The solution of the uncontrolled system~\eqref{eq:EP-uwv} clearly shows that the equilibrium is unstable, whereas that of the controlled system~\eqref{eq:ControlledEP-SE3} stays close to the equilibrium, indicating that the equilibrium is stabilized.

\section{Asymptotic Stabilization}
\subsection{Asymptotic Stabilization by Dissipative Control}
\label{ssec:asymp_stabilization}
Now we would like to introduce an additional dissipative control to have asymptotic stabilization.

We achieve this, as in \cite{BlChLeMaWo2000,BlLeMa2000,BlChLeMa2001}, by applying an additional control $\mathbf{u}^{\rm d}$ to the controlled system~\eqref{eq:ControlledEP2-SE3}:
\begin{equation}
  \label{eq:ControlledEP3-SE3}
  \begin{split}
    \od{}{t} \parentheses{ \pd{\ellc}{\bOmega} }
    &= \pd{\ellc}{\bOmega} \times \bOmega 
    + \pd{\ellc}{\bv} \times \bv
    + \pd{\ellc}{\bGamma} \times \bGamma,
    \\
    \od{}{t} \parentheses{ \pd{\ellc}{\bv} }
    &= \pd{\ellc}{\bv} \times \bOmega + \mathbf{u}^{\rm d},
    \\
    \dot{\bGamma}
    &= \bGamma \times \bOmega.
  \end{split}
\end{equation}
The system is defined on $\R^{3} \times \R^{3} \times \R^{3}$.
However, since $\norm{\bGamma(t)} = 1$ for any time $t$, we consider the system on
\begin{equation*}
  M \defeq \R^{3} \times \R^{3} \times \mathbb{S}^{2}
\end{equation*}
instead.
We shall restrict functions defined on $\R^{3} \times \R^{3} \times \R^{3}$ to $M$ if necessary, but without change of notation for brevity.
We will also write the state variables as $\zeta = (\bOmega, \bv, \bGamma)$ for short in what follows.
Then we have the following result for asymptotic stabilization:

\begin{theorem}[Asymptotic stabilization]
  \label{thm:asymptotic_stability}
  Let $\zeta_{\rm e} \in M$ be an equilibrium of the uncontrolled system~\eqref{eq:EP-uwv}, and $\mathcal{E}\colon M \to \R$ be the Lyapunov function obtained by the energy--Casimir method, i.e., $\mathcal{E}$ is an invariant of \eqref{eq:ControlledEP2-SE3}, $D\mathcal{E}(\zeta_{\rm e}) = 0$, and $D^{2}\mathcal{E}(\zeta_{\rm e})$ is positive definite.
  Let $\mathbf{f}\colon M \to \R^{3}$ be a smooth function, and consider the controlled system~\eqref{eq:ControlledEP3-SE3} with the feedback $\mathbf{u}^{\rm d} = \mathbf{f}(\zeta)$:
  \begin{equation}
    \label{eq:ControlledEP4-SE3}
    \begin{split}
      \od{}{t} \parentheses{ \pd{\ellc}{\bOmega} }
      &= \pd{\ellc}{\bOmega} \times \bOmega 
      + \pd{\ellc}{\bv} \times \bv
      + \pd{\ellc}{\bGamma} \times \bGamma,
      \\
      \od{}{t} \parentheses{ \pd{\ellc}{\bv} }
      &= \pd{\ellc}{\bv} \times \bOmega + \mathbf{f}(\bOmega, \bv, \bGamma),
      \\
      \dot{\bGamma}
      &= \bGamma \times \bOmega,
    \end{split}
  \end{equation}
  and suppose that $\mathbf{f}$ satisfies the following:
  \begin{enumerate}
  \item $\mathbf{f}(\zeta_{\rm e}) = \mathbf{0}$.
    \label{cond:f-1}
    \smallskip
  \item The directional (Lie) derivative $\dot{\mathcal{E}}$ of $\mathcal{E}$ along the solutions of \eqref{eq:ControlledEP4-SE3} gives $\dot{\mathcal{E}}(\zeta) \le 0$ for any $\zeta \in M$.
    \label{cond:f-2}
  \end{enumerate}
  Let $\dot{\mathcal{E}}^{-1}(0) \defeq \setdef{ \zeta \in M }{ \dot{\mathcal{E}}(\zeta) = 0 }$, and $U$ be an open neighborhood of $\zeta_{\rm e}$ such that the only invariant set of \eqref{eq:ControlledEP4-SE3} in $U \cap \dot{\mathcal{E}}^{-1}(0)$ is $U \cap \mathcal{I}$ for some $\mathcal{I} \subset M$.
  Then there exists a compact neighborhood $\Sigma \subset U$ of $\zeta_{\rm e}$ such that any solution to \eqref{eq:ControlledEP4-SE3} starting in $\Sigma$ at $t = 0$ approaches $\Sigma \cap \mathcal{I}$ (which contains $\zeta_{\rm e}$) as $t \to \infty$.
\end{theorem}

\begin{proof}
  Let us first show that $\zeta_{\rm e}$ is an equilibrium of the dissipative controlled system~\eqref{eq:ControlledEP4-SE3}.
  Recall from \cref{ssec:energy-Casimir} that $\zeta_{\rm e}$ is an equilibrium of \eqref{eq:ControlledEP2-SE3} or equivalently \eqref{eq:ControlledEP3-SE3} with $\mathbf{u}^{\rm d} = \mathbf{0}$.
  However, since $\mathbf{f}(\zeta_{\rm e}) = \mathbf{0}$ by assumption, $\zeta_{\rm e}$ is an equilibrium of \eqref{eq:ControlledEP4-SE3} as well.

  Note also that $\mathcal{E}$ is a Lyapunov function for \eqref{eq:ControlledEP4-SE3} as well because the only change due to the dissipative control is that we now have $\dot{\mathcal{E}}(\zeta) \le 0$ instead of $\dot{\mathcal{E}}(\zeta) = 0$.
  So it still implies the Lyapunov stability of $\zeta_{\rm e}$ and hence the existence of a compact neighborhood $\Sigma \subset U$ such that any solution to \eqref{eq:ControlledEP4-SE3} starting in $\Sigma$ at $t = 0$ stays in $\Sigma$ for any $t \ge 0$.
  
  Therefore, LaSalle's Invariance Principle~\cite{La1960} along with the assumption on the invariant set $\mathcal{I}$ implies that any solution starting in $\Sigma$ at $t = 0$ approaches $\Sigma \cap \mathcal{I}$ as $t \to \infty$.

  Note that $\zeta_{\rm e} \in \dot{\mathcal{E}}^{-1}(0)$ due to the condition \ref{cond:f-1}, and so clearly $\zeta_{\rm e} \in \mathcal{I}$ because it is an equilibrium of \eqref{eq:ControlledEP4-SE3}.
  Hence $\zeta_{\rm e} \in \Sigma \cap \mathcal{I}$.
\end{proof}

\subsection{Asymptotic Stabilization of Heavy Top on Movable Base}
Let
\begin{equation}
  \label{eq:mathcalZ-htmb}
  \mathcal{Z}_{\rm e}^{\rm htmb}
  \defeq \setdef{ ( \Omega_{0} \mathbf{E}_{3}, \mathbf{0}, \mathbf{E}_{3} ) }{ \Omega_{0} \in \R }
\end{equation}
be the set of equilibria of the form~\eqref{eq:equilibrium-htmb}.
Notice that each point in this set corresponds to the top spinning at angular velocity $\Omega_{0}$ in the upright position on a stationary base.
  
We shall prove that the solution starting near $\mathcal{Z}_{\rm e}^{\rm htmb}$ at $t = 0$ converges to the point in $\mathcal{Z}_{\rm e}^{\rm htmb}$ determined by setting $\Omega_{0}$ equal to the initial value of $\Omega_{3}$; see \cref{fig:asymp_stability-htmb} below.
This is not quite the asymptotic stability in the conventional sense where any point in a neighborhood of a \textit{single} equilibrium converges to \textit{that} equilibrium.
As we shall explain in \cref{rem:equilibria} below, this subtlety is \textit{not} a drawback of our control law, but is rather due to a nature of this particular control system.
In fact, such a subtlety is not present in the two other examples to follow in the next subsections.

\begin{proposition}[Asymptotic stabilization of heavy top on a movable base]
  \label{prop:asymp_stabilization-htmb}
  Consider the controlled system~\eqref{eq:ControlledEP3-SE3} for the heavy top on a movable base from \cref{ex:htmb}, where
  \begin{equation}
    \label{eq:u^d-htmb}
    \mathbf{u}^{\rm d}
    = \mathcal{N} \left(
      \bv
      + \frac{c\,m l I_1 }{I_1\varrho - m^2l^2} (\bchi \times \bOmega)
    \right) \eqdef \mathbf{f}_{\rm htmb}(\zeta).
  \end{equation}
  with an arbitrary negative-definite matrix $\mathcal{N}$; this is equivalent to applying to the second equation of \eqref{eq:EP-htmb} the control $\mathbf{u} = \mathbf{u}^{\rm p} + \mathbf{u}^{\rm k} + \mathbf{u}^{\rm d}$, where $\mathbf{u}^{\rm p}$ and $\mathbf{u}^{\rm k}$ are those from \cref{prop:stabilization-htmb}.
  For each $\zeta_{\rm e} \in \mathcal{Z}_{\rm e}^{\rm htmb}$, there exists a compact neighborhood $\Sigma \subset \R^{3} \times \R^{3} \times \mathbb{S}^{2}$ of $\zeta_{\rm e}$ such that any solution starting in $\Sigma$ with $\Omega_{3}(0) = \Omega_{0}$ at $t = 0$ approaches the equilibrium $( \Omega_{0} \mathbf{E}_{3}, \mathbf{0}, \mathbf{E}_{3} ) \in \mathcal{Z}_{\rm e}^{\rm htmb}$ as $t \to \infty$.
\end{proposition}
\begin{figure}[htbp]
  \centering
  \includegraphics[width=.375\linewidth]{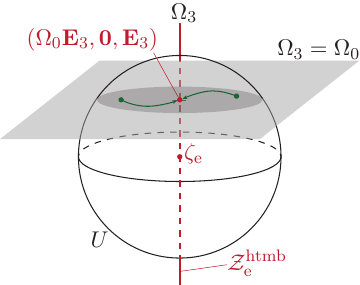}
  \caption{Schematic of \cref{prop:asymp_stabilization-htmb}.
    Subset $\mathcal{Z}_{\rm e}^{\rm htmb}$ is the collection of all upright equilibria of the top spinning with varying value of angular velocity $\Omega_{3}$ on stationary base.
    For any $\zeta_{\rm e} \in \mathcal{Z}_{\rm e}^{\rm htmb}$, one can find a neighborhood $U$ (ball in the figure) so that the following is satisfied:
    For any point in $U$ whose $\Omega_{3}$-value is $\Omega_{0}$, the solution starting at the point converges to equilibrium $( \Omega_{0} \mathbf{E}_{3}, \mathbf{0}, \mathbf{E}_{3} )$ (upright spinning with angular velocity $\Omega_{0}$ on stationary base) as $t \to \infty$.
  }
  \label{fig:asymp_stability-htmb}
\end{figure}
\begin{remark}
  \label{rem:equilibria}
  Notice that $\zeta_{\rm e}$ may not be the same as $( \Omega_{0} \mathbf{E}_{3}, \mathbf{0}, \mathbf{E}_{3} )$.
  The reason for this subtlety is that $\Omega_{3}$ is an invariant of the system even with controls.
  So if $\Omega_{3}(0) \neq \Omega_{0}$, then $\Omega_{3}(t)$ would not converge to $\Omega_{0}$ as $t \to \infty$.
  In other words, the equilibrium to converge to is determined by the initial value of $\Omega_{3}$ as shown in \cref{fig:asymp_stability-htmb}.
  We also emphasize that \textit{we would have the same issue no matter what control $\mathbf{u}$ one applies to the second equation of \eqref{eq:EP-htmb}, because it still gives $\dot{\Omega}_{3} = 0$.
    In other words, one just cannot control the spinning velocity $\Omega_{3}$ of the top however hard one pushes the base.
    So this is rather a nature of this particular control system than an issue specific to our control law.}

  It is an interesting future work to look into the controllability of mechanical systems with broken symmetry in conjunction with the stabilizability discussed here; see \citet{WeBuBu2021} on the controllability of an aerial manipulator as an example of a mechanical system with broken symmetry.
\end{remark}
\begin{proof}[Proof of \cref{prop:asymp_stabilization-htmb}]
  Recall that our control Lyapunov function $\mathcal{E}$ was given in \eqref{eq:mathcalE-htmb}.
  Taking the directional derivative (denoted by $\dot{(\,\cdot\,)}$) of $\mathcal{E}$ along the vector field of the system~\eqref{eq:ControlledEP3-SE3},
  \begin{equation*}
    \dot{\mathcal{E}} = \dot{E}_{\tau, \sigma, \rho} + c\,\dot{E}^{0} + \dot{\Phi} + \dot{\phi},
  \end{equation*}
  and it is easy to see that $\dot{\phi} = \phi'(\Omega_{3}) \dot{\Omega}_{3} = 0$.
  Also, straightforward calculations yield
  \begin{equation*}
    \dot {E}_{\tau, \sigma, \rho} = \bv \cdot \mathbf{u}^{\rm d},
    \qquad
    \dot{E}^{0} = \frac{I_1 m l }{I_1\varrho - m^2 l^2} \bigl(\bchi\times  \bOmega \bigr)\cdot \mathbf{u}^{\rm d}.
  \end{equation*}
  We also have $\dot{\Phi} = 0$ because we have $D_{1}\Phi = D_{2}\Phi = 0$ (see \eqref{eq:Phi-htmb}) as well as $\dot{C}_{3} = 0$.
  Hence we obtain
  \begin{equation}
    \label{eq:dotmathcalE-htmb}
    \dot{\mathcal{E}}
    = \left(
      \bv + \frac{c\,m l I_1 }{I_1\varrho - m^2l^2} (\bchi \times \bOmega) 
    \right) 
    \cdot \mathbf{u}^{\rm d}.
  \end{equation}
  
  Let us consider the feedback control $\mathbf{u}^{\rm d} = \mathbf{f}_{\rm htmb}(\zeta)$ as shown in \eqref{eq:u^d-htmb}.
  Then $\mathbf{f}_{\rm htmb}$ clearly satisfies the conditions~\ref{cond:f-1} and \ref{cond:f-2} on $\mathbf{f}$ stated in \cref{thm:asymptotic_stability}.
  Additionally, \cref{lem:invset-htmb} from \cref{sec:invset} says that there exists a neighborhood $U$ of $\zeta_{\rm e}$ such that $U \cap \mathcal{Z}_{\rm e}^{\rm htmb}$ is the only invariant set in $U \cap \dot{\mathcal{E}}^{-1}(0)$.
  
  Therefore, taking $\mathcal{I} = \mathcal{Z}_{\rm e}^{\rm htmb}$, \cref{thm:asymptotic_stability} implies that there exists a compact neighborhood $\Sigma$ of $\zeta_{\rm e}$ such that any solution starting in $\Sigma$ at $t = 0$ approaches $\Sigma \cap \mathcal{Z}_{\rm e}^{\rm htmb}$ as $t \to \infty$.
  However, since $\Omega_{3}$ is an invariant of the system, this implies that any solution starting in $\Sigma$ with $\Omega_{3}(0) = \Omega_{0}$ approaches the equilibrium $( \Omega_{0} \mathbf{E}_{3}, \mathbf{0}, \mathbf{E}_{3} )$.
\end{proof}

\begin{figure}[htbp]
  \centering
  \subfloat[Body angular velocity $\bOmega$]{
    \includegraphics[width=.31\linewidth]{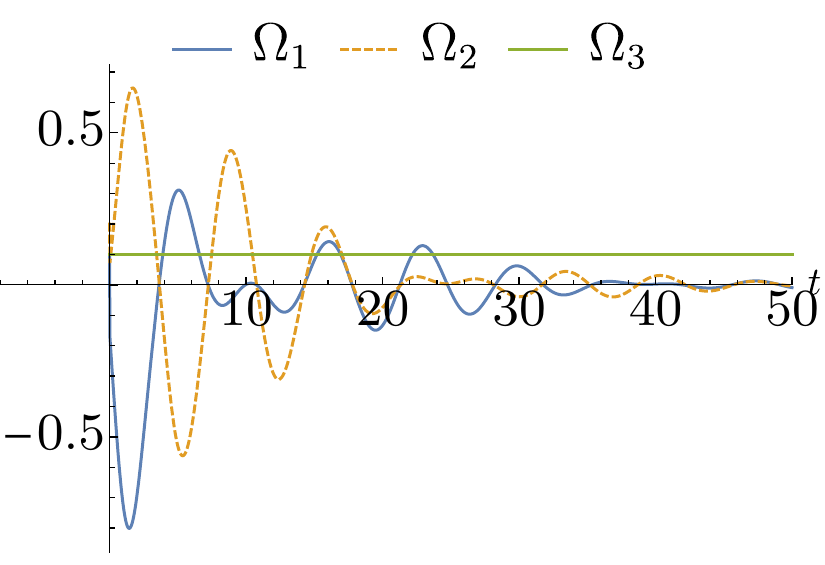}
  }
  \
  \subfloat[Base velocity $\bv$ in body frame]{
    \includegraphics[width=.31\linewidth]{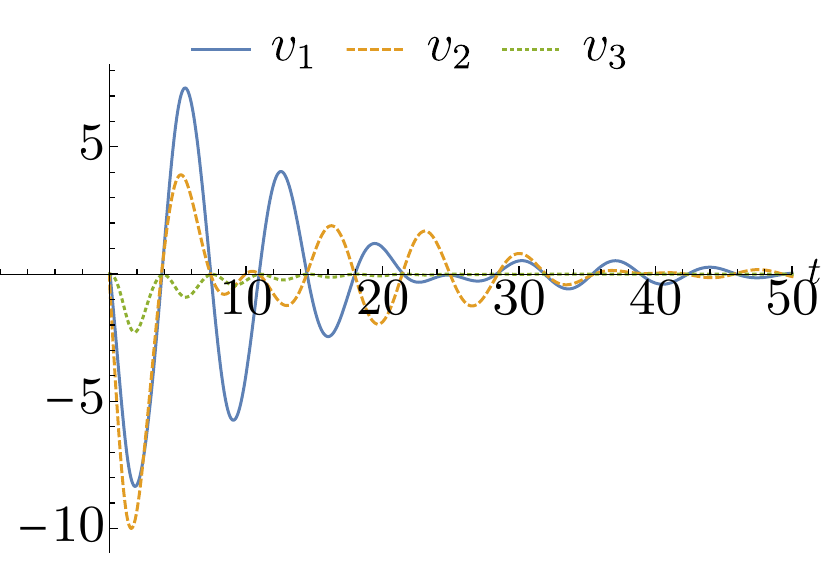}
  }
  \
  \subfloat[Vertical upward direction $\boldsymbol{\Gamma}$ seen from body frame]{
    \includegraphics[width=.31\linewidth]{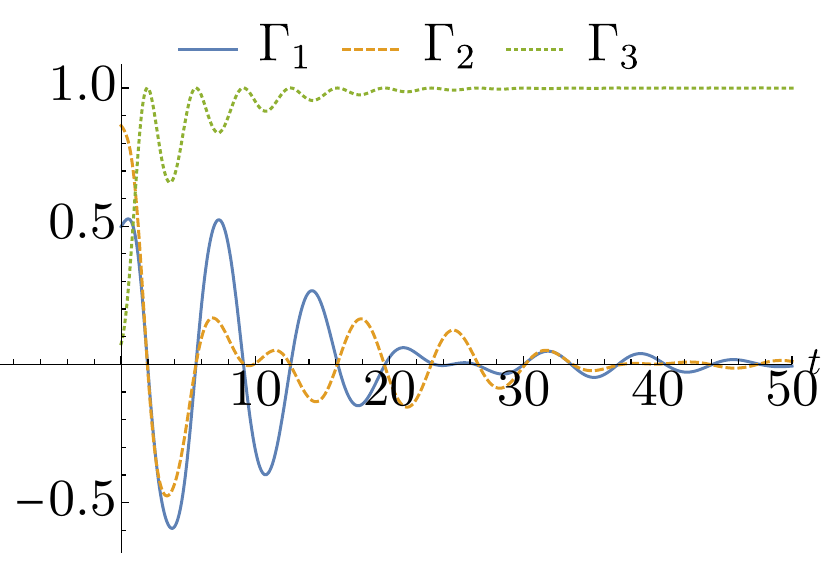}
  }
  \caption{Simulation results for the heavy top on a movable base with dissipative control; the parameters are same as \cref{fig:simulation-HT} except that $\varphi_{0} = \frac{10}{21}\pi$ (near horizontal), $c = 1$, and $\mathcal{N} = -\frac{1}{2}I$.
  The equilibrium is now asymptotically stable.}
  \label{fig:simulation-HT-dis}
\end{figure}

\Cref{fig:simulation-HT-dis} shows the simulation results with the dissipative control, now with $\varphi_{0} = \frac{10}{21}\pi$, i.e., the axis of the top is near the horizontal position.
We see that the control manages to steer the system towards the upright position.

\subsection{Asymptotic Stabilization of Underwater Vehicle}
For the underwater vehicle, we have the asymptotic stability in the conventional sense:
\begin{proposition}[Asymptotic tabilization of underwater vehicle]
  \label{prop:asymp_stabilization-uwv}
  Consider the controlled system~\eqref{eq:ControlledEP3-SE3} for the underwater vehicle from \cref{ex:uwv} where
  \begin{equation}
    \label{eq:u^d-uwv}
    \mathbf{u}^{\rm d}
    = \mathcal{N} \parentheses{ \bv + 2 D_{1}\Phi(C_{1}, C_{2}, C_{3}) \pd{\ellc}{\bv} }
    \eqdef \mathbf{f}_{\rm uwv}(\zeta).
  \end{equation}
  with an arbitrary negative-definite matrix $\mathcal{N}$; this is equivalent to applying to the second equation of \eqref{eq:EP-uwv} the control $\mathbf{u} = \mathbf{u}^{\rm k} + \mathbf{u}^{\rm d}$, where $\mathbf{u}^{\rm k}$ is from \cref{prop:stabilization-uwv}.
  Let
  \begin{equation}
    \label{eq:mathcalZ-uwv}
    \mathcal{Z}_{\rm e}^{\rm uwv}
    \defeq \setdef{ (\mathbf{0}, v_{0}\mathbf{E}_{2}, \mathbf{E}_{3}) }{ v_{0} \in \R\backslash\{0\} }
  \end{equation}
  be the set of equilibria of the form~\eqref{eq:equilibrium-uwv}.
  For each $\zeta_{\rm e} \in \mathcal{Z}_{\rm e}^{\rm uwv}$, use $\Phi$ from \eqref{eq:Phi-uwv} with $C_{1}|_{\rm e} = C_{1}(\zeta_{\rm e})$ and $C_{3}|_{\rm e} = C_{3}(\zeta_{\rm e})$ in the control~\eqref{eq:u^d-uwv}.
  Then, there exists a compact neighborhood $\Sigma \subset \R^{3} \times \R^{3} \times \mathbb{S}^{2}$ of $\zeta_{\rm e}$ such that any solution starting in $\Sigma$ approaches the equilibrium $\zeta_{\rm e}$ as $t \to \infty$.
\end{proposition}
\begin{proof}
  Using the control Lyapunov function~\eqref{eq:mathcalE-uwv}, we have
  \begin{equation}
    \label{eq:dotmathcalE-uwv}
    \dot{\mathcal{E}} = \parentheses{
      \bv + 2 D_{1}\Phi(C_{1}, C_{2}, C_{3}) \pd{\ellc}{\bv}
    } \cdot \mathbf{u}^{\rm d}.
  \end{equation}
  Hence we consider the feedback control $\mathbf{u}^{\rm d} = \mathbf{f}_{\rm uwv}(\zeta)$ as shown in \eqref{eq:u^d-uwv}.
  Then one easily sees that $\mathbf{f}_{\rm uwv}$ satisfies the condition \ref{cond:f-1} on $\mathbf{f}$ stated in \cref{thm:asymptotic_stability} using an expression from \eqref{eq:first_variation-uwv} and $\rho = \diag(m_{1} - \mathcal{K}, m_{2} - \mathcal{K}, m_{3} - \mathcal{K})$.
  It also clearly satisfies the other condition \ref{cond:f-2} by construction.
  
  The rest of the argument is essentially the same as the proof of \cref{prop:asymp_stabilization-htmb} using \cref{lem:invset-uwv}; note however that \cref{lem:invset-uwv} says that the invariant set $\mathcal{I}$ is the equilibrium $\zeta_{\rm e}$ itself as opposed to a family of equilibria.
  Hence \cref{thm:asymptotic_stability} with $\mathcal{I} = \{ \zeta_{\rm e} \}$ gives the desired result.
\end{proof}

\begin{figure}[htbp]
    \includegraphics[width=.32\linewidth]{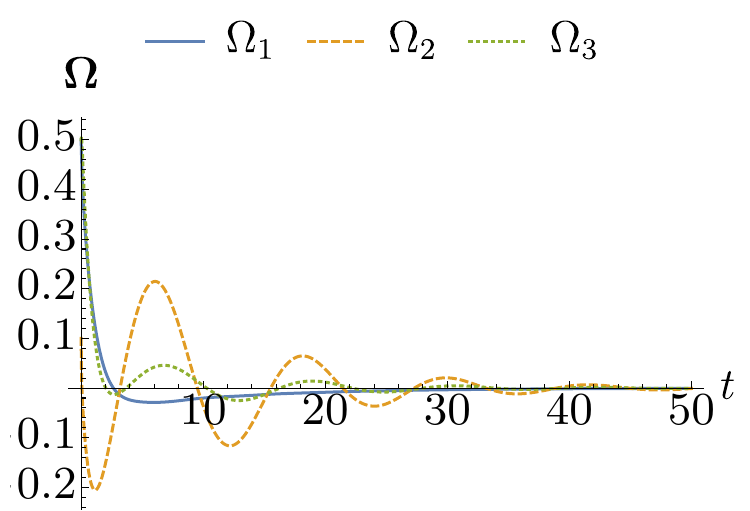}
  \
    \includegraphics[width=.32\linewidth]{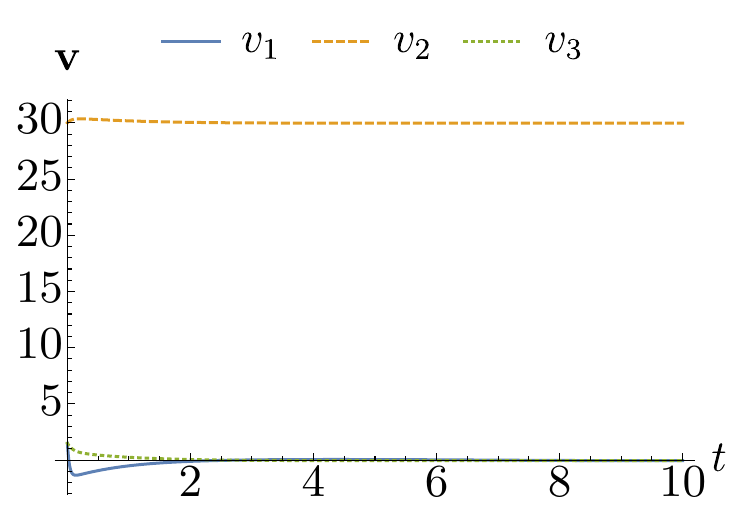}
  \
    \includegraphics[width=.32\linewidth]{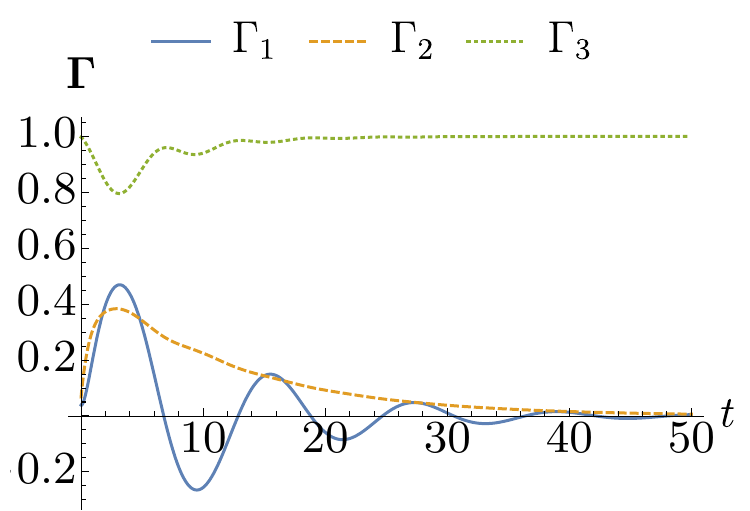}
  \caption{Time evolution of underwater vehicle with dissipative control with $\mathcal{N} = -\diag(2,1,2) \times 10^{6}$.
    Note that $\mathbf{v}$ is plotted for a shorter time interval because it is damped much faster than the other two quantities.}
  \label{fig:simulation-UWV-dis}
\end{figure}

\subsection{Swinging Up the Spherical Pendulum}
As an application of the same control law to a problem with a slightly different flavor, let us consider the problem of swinging up a spherical pendulum on a movable base; see \cref{fig:spmb}.

\begin{figure}
  \centering
  \begin{tikzpicture}[scale=.85]
    \draw [thick, color=gray] (0,1.5) -- (1,5);
    \node[below right] at (1,5) {$m$};
    \node at (.8,3.3) {$l$};
    \shade[ball color=orange] (1,5) circle (5pt); 

    \draw[fill] (0,1.5) circle [radius=0.05];
    \node [draw, cylinder, shape aspect=1, rotate=90, minimum height=.3cm, minimum width=2cm] (c1)  at (0,1.25){};
    \node at (1,0.8)    {$M$};

    \draw[semithick,->]  (-3.5,1) -- (-3.8,0.7) node[below] {$\mathbf{e}_1$};
    \draw[semithick,->]  (-3.5,1) -- (-3,1) node[below right] {$\mathbf{e}_2$};
    \draw[semithick,->]  (-3.5,1) -- (-3.5,1.6) node[above] {$\mathbf{e}_3$};

    \draw[semithick,green!60!black][->]  (-3.5,1) -- (0,1.5);
    \node at (-1.75,1) {\color{green!60!black}{$\mathbf{x}$}};

    \draw[semithick,->]  (0,1.5) -- (-0.525,1.375) node[above] {$\mathbf{E}_1\!$};
    \draw[semithick,->]  (0,1.5) -- (0.375,1.2) node[right] {$\mathbf{E}_2$};
    \draw[semithick,->]  (0,1.5) -- (0.15,2.025) node[right] {$\mathbf{E}_3$};

    \draw[very thick, ->, red!75!black] (-0.75,0.9) -- (-0.25,1.2);
    \node[below right, red!75!black] at (-0.75,1) {$\mathbf{u}$};
  \end{tikzpicture}
  \caption{Spherical pendulum on a movable base. We would like to swing up the pendulum from the (almost) vertical downward position to the upright position by applying a control to the base.}
  \label{fig:spmb}
\end{figure}

Following \cite{ZeLeBl2012}, we treat the pendulum as a degenerate top that does not rotate about its rod.
Specifically, we set the third components of the inertia tensor $\mathbb{I}$ and of the angular velocity $\Omega$ to zero, i.e., $I_3=0$ and $\Omega_3 = 0$.
Assuming that the rod is massless and denoting the bob mass by $m$ and the pendulum length by $l$, the inertia tensor $\mathbb{I}$ becomes $\mathbb{I} =
\bigl[\begin{smallmatrix}
I_1 & 0  \\
0 & I_1 
\end{smallmatrix}
\bigr]
= \bigl[
\begin{smallmatrix}
m l^2 & 0\\
0 & m l^2
\end{smallmatrix} \bigr]$ because we got rid of $\Omega_{3}$ from the formulation.

Since this is the special case of the heavy top with $I_{1} = I_{2} = m l^{2}$, we have the same stability condition under this simplification.
Specifically, we can achieve stability for any $\varrho \in (0,m)$.
Furthermore, since $\Omega_{3} = 0$ here, the set of equilibria $\mathcal{Z}_{\rm e}^{\rm htmb}$ from \eqref{eq:mathcalZ-htmb} becomes a single point.
Hence \cref{prop:asymp_stabilization-htmb} applied to this special case implies the asymptotic stability in the conventional sense:
The solution approaches \textit{the} upright equilibrium as $t \to \infty$.

\begin{figure}[htbp]
  \centering
    \includegraphics[width=.32\linewidth]{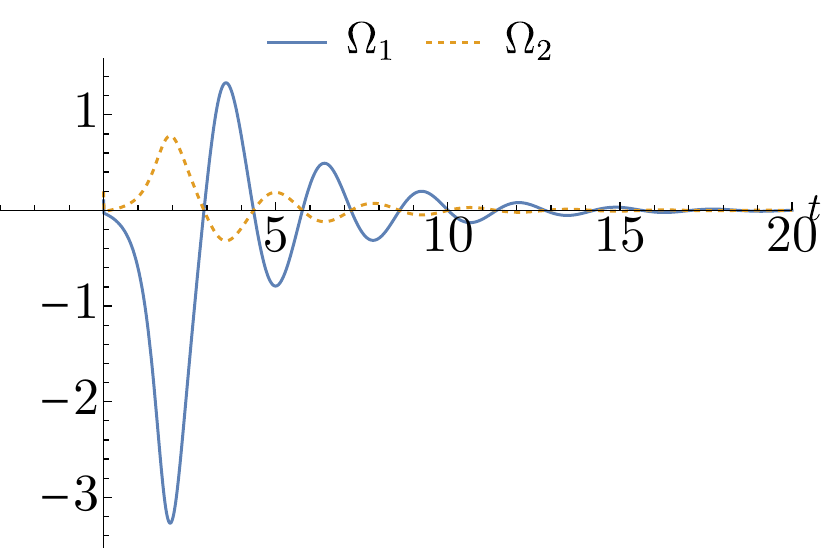}
  \
    \includegraphics[width=.32\linewidth]{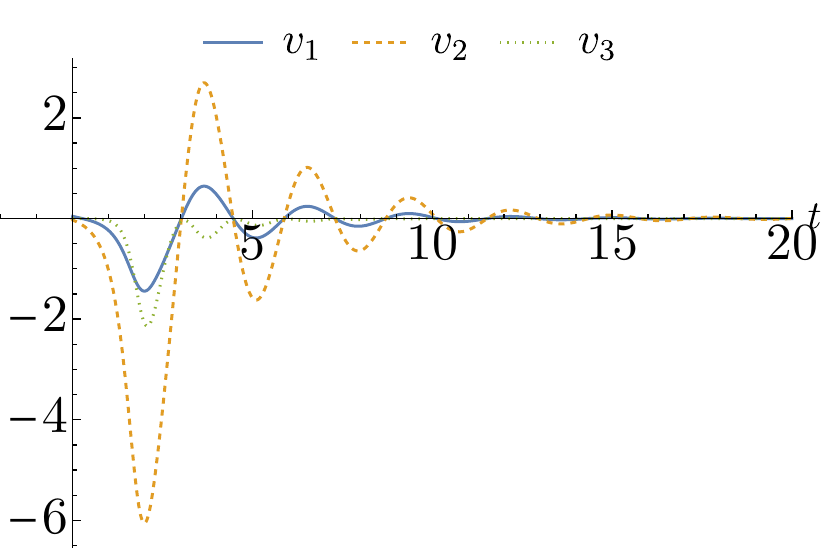}
  \
    \includegraphics[width=.32\linewidth]{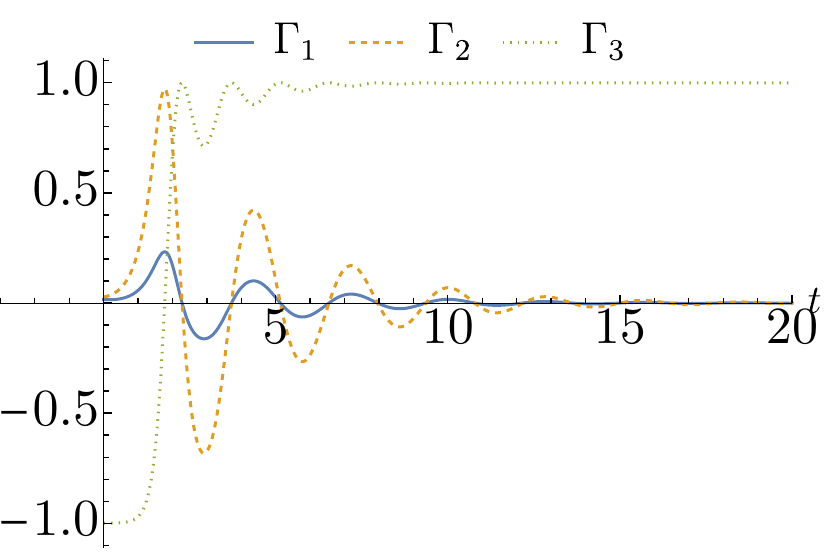}
  \caption{Simulation results for swinging up spherical pendulum on a movable base with $M=0.44\,\mathrm{[kg]}$, $m=0.14\,\mathrm{[kg]}$, $l = 0.215\,\mathrm{[m]}$ (taken from \cite{BlLeMa2000}), and $\varrho = 0.9m$ with initial condition $\bOmega(0) = \mathbf{0}$, $\bv(0) = \mathbf{0}$, and $\boldsymbol{\Gamma}(0) = (\cos\theta_{0} \sin\varphi_{0}, \sin\theta_{0} \sin\varphi_{0}, \cos\varphi_{0})$ with $\theta_{0} = \pi/3$ and $\varphi_{0} = 0.99\pi$, i.e., $\boldsymbol{\Gamma}(0)$ is near the vertical downward position $(0,0,-1)$; also $c = 1$ and $\mathcal{N} = -I$ for dissipative control.
    The pendulum is swung up and asymptotically approaches the upright position $\boldsymbol{\Gamma} = (0,0,1)$.}
  \label{fig:simulation-SP-dis}
\end{figure}

\Cref{fig:simulation-SP-dis} shows the results of simulations.
Note that the initial condition is chosen so that the pendulum is almost downward ($\varphi_{0} = 0.99\pi$) as opposed to exactly downward ($\varphi_{0} = \pi$ or $\boldsymbol{\Gamma}(0) = (0,0,-1)$) because the exact downward position is an equilibrium of the controlled system.
One sees that the pendulum is swung up and asymptotically stabilized towards the upright position.

\section*{Acknowledgments}
We would like to thank Mark Spong for the helpful comments and discussions, Scott Kelly for suggesting us the application to the problem of swinging up the pendulum, and the reviewers for their comments and constructive criticisms.

\appendix
\section{Semidirect Product $\SE(3) \ltimes \R^{4}$}
\label{sec:SE3_ltimes_R4}
This appendix gives a brief summary of the semidirect product Lie groups $\SE(3) \ltimes \R^{4}$ and $\SE(3) \ltimes \R^{3}$ used throughout the paper.

\subsection{$\SE(3)$-action on $\R^{4}$}
\label{ssec:SE3-action_on_R4}
Let $\kappa\colon \SE(3) \to \GL(\R^{4})$ be the left representation of $\SE(3)$ on $\R^{4}$ defined by the standard matrix-vector multiplication:
Writing $s = (R,\mathbf{x}) = \left[
  \begin{smallmatrix}
    R  & \mathbf{x} \\
    \mathbf{0}^{T} & 1
\end{smallmatrix}\right]$,
\begin{equation}
  \label{eq:sigma-R4}
  \kappa(s) y = s y
  = \begin{bmatrix} 
    R  & \mathbf{x} \\
    \mathbf{0}^{T} & 1
  \end{bmatrix}
  \begin{bmatrix}
    \mathbf{y} \\
    \tilde{y}
  \end{bmatrix}
  = \begin{bmatrix}
    R \mathbf{y} + \tilde{y} \mathbf{x}\\ 
    \tilde{y}
  \end{bmatrix}.
\end{equation}
We note in passing that it was also used in the optimal-control formulation of the Kirchhoff elastic rod under gravity by \citet{BoBr2014,BoBr2016}.

Let $(\R^{4})^*$ be the dual of $\R^{4}$.
We identify $(\R^{4})^*$ with $\R^{4}$ via the dot product $\ip{v}{w} \defeq v \cdot w$.
Then the induced left representation $\kappa^{*}\colon \SE(3) \to \GL\parentheses{ (\R^{4})^{*} }$ is defined as
\begin{align*}
  \ip{ \kappa^*(s) \Gamma }{y}
  \defeq \ip{ \kappa(s^{-1})^* \Gamma }{y}
  = \ip{\Gamma}{ \kappa(s^{-1}) y }
  = \ip{\Gamma}{ s^{-1} y }
  = \ip{s^{-T}\Gamma}{y}, 
\end{align*}
and therefore, writing $\Gamma = (\bGamma, h) \in (\R^{4})^{*}$, we have
\begin{equation}
  \label{eq:sigma*-R4}
  \kappa^*(s) \Gamma
  = s^{-T} \Gamma
  = \begin{bmatrix}
    R & 0 \\
    -\mathbf{x}^{T}R & 1
  \end{bmatrix}
  \begin{bmatrix}
    \bGamma \\
    h
  \end{bmatrix}
  =
  \begin{bmatrix}
    R\bGamma \\
    -\mathbf{x}^{T}R \bGamma + h
  \end{bmatrix}.
\end{equation}

We may identify the Lie algebra $\se(3) = \so(3) \ltimes \R^{3}$ with $\R^{3} \times \R^{3}$ via the hat map~\eqref{eq:hat_map}:
\begin{equation*}
  \begin{bmatrix}
    \hat{\Omega} & \bv \\
    \mathbf{0}^{T} & 0
  \end{bmatrix}
  \in \se(3)
  \leftrightarrow
  (\bOmega, \bv) \in \R^{3} \times \R^{3}.
\end{equation*}
Then we may write the induced action of $\se(3)$ on $\R^4$ as
\begin{equation}
  \label{eq:sigma'-R4}
  \kappa'(\bOmega,\bv) y
  = \begin{bmatrix} 
    \hat{\Omega}  & \bv\\
    \mathbf{0}^{T} & 0
  \end{bmatrix}
  \begin{bmatrix}
    \mathbf{y} \\
    \tilde{y}
  \end{bmatrix} 
  =
  \begin{bmatrix}
    \hat{\Omega} \mathbf{y} + \tilde{y} \bv\\
    0
  \end{bmatrix}
  =\begin{bmatrix}
    \bOmega \times \mathbf{y} + \tilde{y} \bv\\
    0
  \end{bmatrix}.
\end{equation}
This induces the Lie algebra action on the dual $(\R^{4})^{*}$ as follows:
\begin{align*}
  \ip{ \kappa'(\bOmega,\bv)^*\Gamma }{y}
  = \ip{\Gamma}{ \kappa'(\bOmega,\bv) y }
  = \mathbf{\Gamma} \cdot (\bOmega \times \mathbf{y} + \tilde{y} \bv)
  = (\mathbf{\bGamma} \times \bOmega) \cdot \mathbf{y} + (\bGamma \cdot \bv) \tilde{y}
\end{align*}
that is,
\begin{equation}
  \label{eq:sigma'star-R4}
  \kappa'(\bOmega,\bv)^* \Gamma =
  \begin{bmatrix}
    \bGamma \times \bOmega \\
    \bGamma \cdot \bv
  \end{bmatrix}.
\end{equation}

For any $y \in \R^{4}$, define a linear map $\kappa_{y} \colon \se(3) \to \R^{4}$ by
\begin{equation*}
  \kappa_{y}(\bOmega,\bv) \defeq \kappa'(\bOmega,\bv) y =
  \begin{bmatrix}
    \bOmega \times \mathbf{y} + \tilde{y} \bv\\
    0
  \end{bmatrix}.
\end{equation*}
Using its dual $\kappa_{y}^{*}\colon (\R^{4})^{*} \to \se(3)^*$, we define the momentum map $\mathbf{K}\colon \R^{4} \times (\R^{4})^{*} \to \se(3)^{*}$ as 
\begin{equation*}
  \mathbf{K}(y, \Gamma) \defeq \kappa^{*}_{y}(\Gamma).
\end{equation*}
To make it more concrete, let us identify $\se(3)^{*}$ with $\se(3)$ via the following inner product on $\se(3)$:
\begin{equation*}
  \ip{ (\bOmega,\bv) }{ (\bXi,\mathbf{w}) }
  \defeq
  \bOmega \cdot \boldsymbol{\Xi}
  + \bv \cdot \mathbf{w}.
\end{equation*}
Then we have, for any $(\bOmega,\bv) \in \se(3)$ and $(y, \Gamma) \in \R^{4} \times (\R^{4})^{*}$, 
\begin{align*}
  \ip{ \mathbf{K}(y, \Gamma) }{ (\bOmega,\bv) }
  = \ip{ \kappa^{*}_{y}(\Gamma) }{ (\bOmega,\bv) } 
  &= \ip{\Gamma}{ \kappa_{y}(\bOmega,\bv) } \\
  &= \bGamma \cdot (\bOmega \times \mathbf{y} + \tilde{y}\bv)
  = (\mathbf{y} \times \bGamma) \cdot \bOmega + (\tilde{y}\bGamma) \cdot \bv,
\end{align*}
which gives
\begin{equation}
  \label{eq:mathbfK-R4}
  \mathbf{K}(y,\Gamma) = ( \mathbf{y} \times \bGamma, \tilde{y}\bGamma).
\end{equation}

\subsection{$\SE(3)$-action on $\R^{3}$}
\label{ssec:SE3-action_on_R3}
Setting $\tilde{y} = 0$ above yields the representation
\begin{equation}
  \label{eq:sigma-R3}
  \kappa\colon \SE(3) \to \GL(\R^{3});
  \qquad
  \kappa(s) \mathbf{y} = \kappa(R,\mathbf{x}) \mathbf{y} = R\mathbf{y}.
\end{equation}
Note that this is \textit{not} the standard $\SE(3)$-action on $\R^{3}$ by rotation and translation.
As a result, we have
\begin{equation}
  \label{eq:sigmas-R3}
  \begin{array}{c}
    \DS \kappa^{*}(R,\mathbf{x}) \bGamma = R\bGamma,
    \qquad
    \kappa'(\bOmega,\bv) \mathbf{y}
    = \kappa_{\mathbf{y}}(\bOmega,\bv)
    = \bOmega \times \mathbf{y},
    \medskip\\
    \DS \kappa'(\bOmega,\bv)^* \bGamma = \bGamma \times \bOmega,
    \qquad
    \mathbf{K}(\mathbf{y},\bGamma) = \mathbf{y} \times \bGamma.
  \end{array}
\end{equation}

\subsection{Lie brackets and coadjoint operator}
Let us find the Lie bracket associated with the semidirect product Lie algebra $\se(3)\ltimes \R^4$. 
Let $(\zeta_1, w_1), (\zeta_2, w_2) \in \se(3) \ltimes \R^4$, where
\begin{equation*}
  (\zeta_1, w_1)=\big( (\bOmega_1,\bv_1), (\mathbf{w}_1,\tilde{w}_1) \big),
  \qquad
  (\zeta_2, w_2)=\big( (\bOmega_2,\bv_2), (\mathbf{w}_2,\tilde{w}_2) \big).
\end{equation*}
Then, the Lie bracket is given by, using $\kappa'$ from \eqref{eq:sigma'-R4} (see also \eqref{eq:ad-sdp}),
\begin{align*}
  [ (\zeta_1, w_1), (\zeta_2, w_2)] 
  &= \ad_{(\zeta_{1},w_{1})}(\zeta_{2},w_{2})
  = \Big(
  [\zeta_1,\zeta_2], \kappa'(\zeta_1) w_2 - \kappa'(\zeta_2) w_1
  \Big)
  \nonumber\\
  &= \Big(
    \big [
    (\bOmega_1,\bv_1),  (\bOmega_2,\bv_2)
    \big ],
    \kappa'(\bOmega_1,\bv_1) (\mathbf{w}_2,\tilde{w}_2)
  - \kappa'(\bOmega_2,\bv_2) (\mathbf{w}_1,\tilde{w}_1)
  \Big)
  \nonumber\\
  &= \Big(
    \big (
    \bOmega_1 \times \bOmega_2,
    \bOmega_1 \times \bv_2 - \bOmega_2 \times \bv_1
    \big), \nonumber\\
  &\qquad
    (\bOmega_1 \times \mathbf{w}_2  + \tilde{w}_2 \bv_1
    - \bOmega_2 \times \mathbf{w}_1 - \tilde{w}_1 \bv_2,
    0 )
    \Big),
\end{align*}
and for $\se(3) \ltimes \R^{3}$,
\begin{equation*}
  \Big [
  \big( 
  \zeta_{1}, \mathbf{w}_1
  \big ),
  \big(
  \zeta_{2}, \mathbf{w}_2
  \big)
  \Big ]
  =
  \Big(
  \big (
  \bOmega_1 \times \bOmega_2,
  \bOmega_1 \times \bv_2 - \bOmega_2 \times \bv_1
  \big), \\
  \bOmega_1 \times \mathbf{w}_2
  - \bOmega_2 \times \mathbf{w}_1 
  \Big).
\end{equation*}

\subsection{Lie--Poisson brackets on $(\se(3) \ltimes \R^{4})^{*}$ and $(\se(3) \ltimes \R^{3})^{*}$}
\label{ssec:LPB}
Using coordinates
\begin{equation*}
  (\mu, \Gamma) = ((\bPi,\bP), (\bGamma,h))
\end{equation*}
for $(\se(3) \ltimes \R^{4})^{*} \cong \R^{3} \times \R^{3} \times \R^{4}$, the $(-)$-Lie--Poisson bracket on $(\se(3) \ltimes \R^{4})^{*}$ is given by
\begin{equation*}
  \PB{F}{G} (\mu,\Gamma)
  = -\ip{ (\mu,\Gamma) }{ \brackets{ \fd{F}{(\mu,\Gamma)}, \fd{G}{(\mu,\Gamma)} } },
\end{equation*}
where, using the Lie brackets from the previous subsection,
\begin{align*}
  \brackets{ \fd{F}{(\mu,\Gamma)}, \fd{G}{(\mu,\Gamma)} }
  &= \Biggl(
    \parentheses{
    \pd{F}{\bPi} \times \pd{G}{\bPi},\,
    \pd{F}{\bPi} \times \pd{G}{\bP} - \pd{G}{\bPi} \times \pd{F}{\bP}
    },\\
  &\qquad
    \parentheses{
    \pd{F}{\bPi} \times \pd{G}{\bGamma}  + \pd{G}{h} \pd{F}{\bP}
    - \pd{G}{\bPi} \times \pd{F}{\bGamma} - \pd{F}{h} \pd{G}{\bP},\,
    0
    }
    \Biggr).
\end{align*}
Hence more concretely,
\begin{align*}
  \PB{F}{G} ( (\bPi, \bP), (\bGamma,h) )
  &= -\bPi \cdot \parentheses{
    \pd{F}{\bPi} \times \pd{G}{\bPi}
    }
    - \bP \cdot \parentheses{
    \pd{F}{\bPi} \times \pd{G}{\bP} - \pd{G}{\bPi} \times \pd{F}{\bP}
    }
  \\
  &\quad - \bGamma \cdot \parentheses{
    \pd{F}{\bPi} \times \pd{G}{\bGamma} + \pd{G}{h} \pd{F}{\bP}
    - \pd{G}{\bPi} \times \pd{F}{\bGamma} - \pd{F}{h} \pd{G}{\bP}
    }.
\end{align*}
Then, one easily sees that $C_1=\|\boldsymbol \Gamma \|^2$ and $C_2=\|\bP \times \boldsymbol \Gamma \|^2$ are Casimirs, i.e., $\PB{F}{C_{i}} = 0$ for any $F \in C^{\infty}\big( (\se(3) \ltimes \R^{4})^{*} \big)$ and $i = 1, 2$.

On the other hand, the Lie--Poisson bracket on $(\se(3) \ltimes \R^{3})^{*}$ is 
\begin{align*}
  \{F,G \} ( (\bPi, \bP), \bGamma )
  &= -\bPi \cdot \parentheses{
    \pd{F}{\bPi} \times \pd{G}{\bPi}
    }
    - \bP \cdot \parentheses{
    \pd{F}{\bPi} \times \pd{G}{\bP} - \pd{G}{\bPi} \times \pd{F}{\bP}
    }
  \\
  &\quad - \bGamma \cdot \parentheses{
    \pd{F}{\bPi} \times \pd{G}{\bGamma}
    - \pd{G}{\bPi} \times \pd{F}{\bGamma}
    }.
\end{align*}
In this case, we have an additional Casimir:
\begin{equation}
  \label{eq:Casimirs-LP}
  C_1=\|\bP \|^2,
  \qquad
  C_2= \bP \cdot \bGamma ,
  \qquad
  C_3=\|\bGamma \|^2.
\end{equation}

\section{Some Lemmas on Invariant Sets}
\label{sec:invset}
\begin{lemma}
  \label{lem:invset-htmb}
  Consider the system~\eqref{eq:ControlledEP3-SE3} with the dissipative control \eqref{eq:u^d-htmb} for the heavy top on a movable base (\cref{ex:htmb}), and define the set
  \begin{equation*}
    \dot{\mathcal{E}}^{-1}(0) \defeq \setdef{ \zeta \in \R^{3} \times \R^{3} \times \mathbb{S}^{2} }{ \dot{\mathcal{E}}(\zeta) = 0 }
  \end{equation*}
  with the Lyapunov function~\eqref{eq:mathcalE-htmb}.
  Then, for each equilibrium $\zeta_{\rm e} \in \mathcal{Z}_{\rm e}^{\rm uwv}$ (defined in \eqref{eq:mathcalZ-htmb}), there exists an open neighborhood $U$ of $\zeta_{\rm e}$ such that the only invariant set inside $U \cap \dot{\mathcal{E}}^{-1}(0)$ is $U \cap \mathcal{Z}_{\rm e}^{\rm htmb}$.
\end{lemma}
\begin{proof}
  In view of \eqref{eq:u^d-htmb} and \eqref{eq:dotmathcalE-htmb}, we have
  \begin{equation*}
    \dot{\mathcal{E}}
    = \left(
      \bv + \frac{c\,m l I_1 }{I_1\varrho - m^2l^2} (\bchi \times \bOmega) 
    \right)^{T}
    \mathcal{N}
    \left(
      \bv + \frac{c\,m l I_1 }{I_1\varrho - m^2l^2} (\bchi \times \bOmega) 
    \right),
  \end{equation*}
  and since $\mathcal{N}$ is assumed to be negative-definite,
  \begin{equation}
    \label{eq:constraint-htmb}
    \dot{\mathcal{E}}(\zeta) = 0
    \iff
    \bv = -\frac{c\,m l I_1 }{I_1\varrho - m^2l^2} (\bchi \times \bOmega)
    \iff \mathbf{f}_{\rm htmb}(\zeta) = \mathbf{0}.
  \end{equation}
  Since $\bchi = (0,0,1)$, this implies that $v_{3} = 0$.
  But then the equations of motion satisfying \eqref{eq:constraint-htmb} gives
  \begin{equation*}
    \dot{v}_{3} = m l \frac{(1+c) I_{1} \varrho - m^{2} l^{2}}{I_{1} \varrho - m^{2} l^{2}} (\Omega_{1}^{2} + \Omega_{2}^{2}).
  \end{equation*}
  The fraction on the right-hand side is non-zero because of the condition on $\varrho$ from \eqref{eq:posdef-htmb}.
  Hence the solution in the invariant set necessarily satisfies $\Omega_{1} = \Omega_{2} = 0$.
  It also implies via \eqref{eq:constraint-htmb} that $v_{1} = v_{2} = 0$ as well, i.e., $\mathbf{v} = \mathbf{0}$.
  Then the equations of motion now give
  \begin{equation*}
    \dot{\Omega}_{1} = m \mathrm{g} l \Gamma_{2},
    \qquad
    \dot{\Omega}_{2} = -m \mathrm{g} l \Gamma_{1}.
  \end{equation*}
  However, since $\Omega_{1} = \Omega_{2} = 0$, we have $\Gamma_{1} = \Gamma_{2} = 0$, and hence $\Gamma_{3} = \pm1$ because $\norm{\bGamma} = 1$.
  We may then take a neighborhood $U$ of $\zeta_{\rm e}$ to exclude $\Gamma_{3} = -1$.
  As a result, $\bGamma = (0,0,1)$, and thus $U \cap \mathcal{Z}_{\rm e}^{\rm htmb}$ is the only invariant set in $U \cap \dot{\mathcal{E}}^{-1}(0)$.
\end{proof}

\begin{lemma}
  \label{lem:invset-uwv}
  Consider the system~\eqref{eq:ControlledEP3-SE3} with the dissipative control \eqref{eq:u^d-uwv} for the underwater vehicle (\cref{ex:uwv}), and define the set $\dot{\mathcal{E}}^{-1}(0)$ with the Lyapunov function $\mathcal{E}$ from \eqref{eq:mathcalE-uwv}.
  Then, for each equilibrium $\zeta_{\rm e} \in \mathcal{Z}_{\rm e}^{\rm uwv}$ (defined in \eqref{eq:mathcalZ-uwv}), there exists an open neighborhood $U$ of $\zeta_{\rm e}$ such that the only invariant set inside $U \cap \dot{\mathcal{E}}^{-1}(0)$ is $\{\zeta_{\rm e}\}$.
\end{lemma}
\begin{proof}
  In view of \eqref{eq:u^d-uwv} and \eqref{eq:dotmathcalE-uwv}, one sees
  \begin{equation*}
    \dot{\mathcal{E}}(\zeta) = 0
    \iff
    \bv = -2 D_{1}\Phi(C_{1}, C_{2}, C_{3}) \pd{\ellc}{\bv}
    \iff \mathbf{f}_{\rm uwv}(\zeta) = \mathbf{0}.
  \end{equation*}
  However, using \eqref{eq:Phi-uwv},
  \begin{equation*}
    D_{1}\Phi(C_{1}, C_{2}, C_{3})
    = \frac{C_{1} - C_{1}|_{\rm e}}{(\mathcal{K} - m_{2})^{3} v_{0}^{2}}
    + \frac{1}{2(\mathcal{K} - m_{2})}.
  \end{equation*}
  Also, using the expression for $C_{1}$ in \eqref{eq:Casimirs} and \eqref{eq:ControlledEP4-SE3}, we see that
  \begin{equation*}
    \dot{C}_{1} = 2 \pd{\ellc}{\bv} \cdot \mathbf{f}_{\rm uwv}(\zeta) = 0,
  \end{equation*}
  because we have $\mathbf{f}_{\rm uwv}(\zeta) = \mathbf{0}$ now.
  Hence $D_{1}\Phi(C_{1}, C_{2}, C_{3})$ is constant if $\dot{\mathcal{E}} = 0$.
  Now, since
  \begin{equation*}
    \pd{\ellc}{\bv} = \rho \mathbf{v} + m l (\bOmega \times \bchi),
  \end{equation*}
  with $\rho = \diag(m_{1} - \mathcal{K}, m_{2} - \mathcal{K}, m_{3} - \mathcal{K})$, we have
  \begin{multline}
    \bv = -2 D_{1}\Phi(C_{1}, C_{2}, C_{3}) \parentheses{ \rho \mathbf{v} + m l (\bOmega \times \bchi) } \\
    \iff
    (I + 2 D_{1}\Phi(C_{1}, C_{2}, C_{3}) \rho ) \mathbf{v} = -2 m l\, D_{1}\Phi(C_{1}, C_{2}, C_{3}) (\bOmega \times \bchi).
    \label{eq:constraint-uwv}
  \end{multline}
  Particularly, since $\bchi = (0, 0, -1)$, the second and third components give
  \begin{equation*}
    \frac{C_{1} - C_{1}|_{\rm e}}{(\mathcal{K} - m_{2})^{2} v_{0}^{2}}\, v_{2}
    = -m l \parentheses{
      \frac{C_{1} - C_{1}|_{\rm e}}{(\mathcal{K} - m_{2})^{3} v_{0}^{2}}
      + \frac{1}{2(\mathcal{K} - m_{2})}
    } \Omega_{1},
    \qquad
    v_{3} = 0.
  \end{equation*}
  
  In what follows, we consider two cases depending on the value of $C_{1}$, which takes the form
  \begin{equation*}
    C_{1} = \parentheses{ (m_{1}-\mathcal{K}) v_{1} - m l \Omega_{2} }^{2}
    + \parentheses{ (m_{2}-\mathcal{K}) v_{2} + m l \Omega_{1} }^{2}
    + \parentheses{ (m_{3}-\mathcal{K}) v_{3} }^{2}.
  \end{equation*}
  
  \begin{enumerate}[leftmargin=.75em, itemindent=3em, label=Case~\arabic*:]
  \item $C_{1} \neq C_{1}|_{\rm e}$\\
    In this case, one may use \eqref{eq:constraint-uwv} to express $\Omega_{1}$ and $\Omega_{2}$ as constant multiples of $v_{2}$ and $v_{1}$, respectively. 
    Then the equations of motion give
    \begin{equation}
      \label{eq:veq}
      \od{}{t}
      \begin{bmatrix}
        v_{1} \\
        v_{2}
      \end{bmatrix}
      = \frac{v_{0}^{2}(\mathcal{K} - m_{2})^{3}}{2C_{1} - (\mathcal{K} - m_{2})^{2}v_{0}^{2}}\, \Omega_{3}
      \begin{bmatrix}
        -v_{2} \\
        v_{1}
      \end{bmatrix}.
    \end{equation}
    and hence
    \begin{equation}
      \label{eq:vsq}
      v_{1}^{2} + v_{2}^{2} = \bar{v}^{2}
      \text{ for some $\bar{v} \in \R$}.
    \end{equation}
    It also gives
    \begin{equation*}
      \dot{v}_{3} = a_{1} v_{1}^{2} + a_{2} v_{2}^{2}
    \end{equation*}
    with some constants $a_{1}, a_{2} \in \R$ with $a_{1} \neq a_{2}$.
    However, because $v_{3} = 0$, we have $a_{1} v_{1}^{2} + a_{2} v_{2}^{2} = 0$; this along with \eqref{eq:vsq} then implies that both $v_{1}$ and $v_{2}$ are constant.
    Therefore, the right-hand side of \eqref{eq:veq} vanishes, i.e., either $\Omega_{3} \neq 0$ and $v_{1} = v_{2} = 0$ or $\Omega_{3} = 0$.
    
    In the former case, we also have $\Omega_{1} = \Omega_{2} = 0$, and this implies that $C_{1} = 0$.
    Therefore, we may take a small enough neighborhood $U$ of the equilibrium (at which $C_{1} = C_{1}|_{\rm e} = (\mathcal{K} - m_{2})^{2} v_{0}^{2} > 0$) so that $C_{1} > 0$ on $U$ to exclude this case.
    
    In the latter case, setting $\Omega_{3} = 0$ gives
    \begin{equation*}
      \dot{\Omega}_{1} = -m \mathrm{g} l \Gamma_{2},
      \qquad
      \dot{\Omega}_{2} = m \mathrm{g} l \Gamma_{1}.
    \end{equation*}
    However, since $\Omega_{1}$ and $\Omega_{2}$ are constant, we have $\Gamma_{1} = \Gamma_{2} = 0$, and so $\Gamma_{3} = \pm1$.
    Either way, setting $\dot{v}_{1} = \dot{v}_{2} = 0$ again leads to $v_{1} = v_{2} = 0$ because of \eqref{eq:veq}, and hence $\Omega_{1} = \Omega_{2} = 0$.
    As a result, $C_{1} = 0$ again and so we may exclude this case as well.
    \smallskip
  \item Case~2: $C_{1} = C_{1}|_{\rm e}$\\
    In this case, \eqref{eq:constraint-uwv} gives
    \begin{equation*}
      \Omega_{1} = 0,
      \qquad
      \Omega_{2} = \frac{m_{1} - m_{2}}{m l}\,v_{1},
      \qquad
      v_{3} = 0,
    \end{equation*}
    and as a result, the equations of motion satisfying \eqref{eq:constraint-uwv} gives
    \begin{equation*}
      \dot{v}_{3} = \frac{ (m_{1} - m_{2})(m_{2} - \mathcal{K}) }{ m l }\,v_{1}^{2}.
    \end{equation*}
    Since $m_{1} \neq m_{2}$ for our case and $\mathcal{K} > m_{2}$ because of \eqref{eq:mathcalK-stability} as well as $v_{3} = 0$, it follows that $v_{1} = 0$, and hence $\Omega_{2} = 0$ as well.
    This in turn implies
    \begin{equation*}
      \dot{\Omega}_{1} = -m \mathrm{g} l \Gamma_{2},
      \qquad
      \dot{\Omega}_{2} = m l (\mathrm{g} \Gamma_{1} - v_{2} \Omega_{3}),
      \qquad
      \dot{v}_{1} = (m_{2} - \mathcal{K}) v_{2} \Omega_{3}.
    \end{equation*}
    The first equation with $\Omega_{1} = 0$ from above implies $\Gamma_{2} = 0$.
    The second with $\Omega_{2} = 0$ implies $v_{2} \Omega_{3} = \mathrm{g} \Gamma_{1}$, and substituting this to the last equation from above,
    \begin{equation*}
      \dot{v}_{1} = (m_{2} - \mathcal{K}) \mathrm{g} \Gamma_{1},
    \end{equation*}
    but then, since $v_{1} = 0$, we have $\Gamma_{1} = 0$, and so $v_{2} \Omega_{3} = 0$ as well.
    Since $\Gamma_{1} = \Gamma_{2} = 0$, we have $\Gamma_{3} = \pm1$.
    Taking a small enough neighborhood $U$ of the equilibrium to exclude $\Gamma_{3} = -1$, we have $\Gamma_{3} = 1$.
    Now, since $v_{1} = v_{3} = \Omega_{1} = \Omega_{2} = 0$, we have $C_{1} = (\mathcal{K} - m_{2})^{2} v_{2}^{2}$.
    However, we are assuming that $C_{1} = C_{1}|_{\rm e} = (\mathcal{K} - m_{2})^{2} v_{0}^{2}$, we have $v_{2} = \pm v_{0} \neq 0$.
    Again, taking $U$ small enough, we have $v_{2} = v_{0}$.
  \end{enumerate}
\end{proof}

\bibliography{CtrlEPwithBSym1}
\bibliographystyle{plainnat}
                                        
\end{document}